\documentclass[10pt,twoside]{article}
\usepackage{amsmath,amssymb}
\usepackage{graphicx}
\usepackage[footnotesize]{subfigure}
\usepackage[footnotesize]{caption2}
\usepackage{authblk}

\allowdisplaybreaks[4]

\newtheorem{lemma}{{\bf Lemma}}[section]

\newtheorem{theorem}{{\bf Theorem}}[section]
\newtheorem{remark}{{\bf Remark}}[section]
\newtheorem{proposition}{{\bf Proposition}}[section]

\makeatletter
\evensidemargin
\oddsidemargin
\makeatother

\title{
\Large\bf Double canard cycles in singularly perturbed planar systems\\ with two canard points
}

\author{{\sc Shuang Chen$^{a,b}$}, {\sc Jinqiao Duan$^{c}$}, {\sc Ji Li$^{a,}$}
\\
{\small $^{a}$ School of Mathematics and Statistics, Huazhong University of Sciences and Technology}\\
{\small Wuhan, Hubei 430074, P. R. China}\\
{\small $^{b}$ Center for Mathematical Sciences, Huazhong University of Sciences and Technology}\\
{\small Wuhan, Hubei 430074, P. R. China}\\
{\small $^{c}$ Department of Applied Mathematics, Illinois Institute of Technology}\\
{\small Chicago, IL 60616, USA}
}

\date{}

\begin{document}
\maketitle
\begin{abstract}

We consider double canard cycles including two canards in singularly perturbed planar  systems with two canard points.
Previous work studied the complex oscillations including relaxation oscillations and canard cycles
in singularly perturbed planar systems with one-parameter layer equations,
which have precisely one canard point, two jump points or one canard point and one jump point.
Based on the normal form theory, blow-up technique and Melnikov theory,
we investigate double canard cycles induced by two Hopf breaking mechanisms at two non-degenerate canard points.
Finally, we apply the obtained results to a class of cubic Li\'{e}nard equations with quadratic damping.


\vskip 0.2cm
{\bf AMS(2020) Subject Classification}: 34C05; 34E17; 34E20.
\end{abstract}
\baselineskip 15pt
\parskip 10pt

\thispagestyle{empty}
\setcounter{page}{1}


\section{Introduction}
\setcounter{equation}{0}
\setcounter{lemma}{0}
\setcounter{remark}{0}

Singular perturbation problems induced by multiple time scales widely appear  in applied science and engineering,
such as cellular physiology, fluid mechanics, population dynamics  and so on \cite{Keener-Sneyd-98,Kevorkan-cole-81,Kuehn-15,Rubin-Terman}.
These singularly perturbed systems are also referred to as slow-fast systems.
Based on the normally hyperbolic invariant manifolds theory (see, for instance, \cite{Fenichel-71,Fenichel-74,Fenichel-77,Wiggins-94}),
Fenichel \cite{Fenichel-79} in 1979 laid the foundation of geometric singular perturbation theory (abbreviated as GSPT)
to investigate multiple time scales dynamics.
Since then, GSPT has become an hotspot research subject in the field of dynamical systems.
There is an enormous literature on this topic.
We refer the readers to  the excellent survey articles
\cite{Deng-03,Deng-17,DuLiLi-18,Dumortieretal-Roussarie-96,Bo-Liu-07,Hek-2010,Jones-95,Li-llibre-12,Schecter2-08} and the references therein.

In this paper we consider a singularly perturbed planar system in the form
\begin{eqnarray}
\label{slow-1}
\begin{split}
\varepsilon\frac{d x}{d \tau} &=  \varepsilon \dot x = f(x,y,\mu,\varepsilon),
\\
\frac{d y}{d \tau} &= \dot y  = g(x,y,\mu,\varepsilon),
\end{split}
\end{eqnarray}
where $(x,y)\in \mathbb{R}^{2}$,
$\mu=(\lambda,\eta)=(\lambda,\eta_{1},...,\eta_{m})\in \mathbb{R}\times \mathbb{R}^{m}$ with $m\geq 1$,
the parameter $\varepsilon$ satisfies $0<\varepsilon\leq \varepsilon_{0}\ll1$ for some small $\varepsilon_{0}$,
and the functions $f$ and $g$ are $C^{k}$ with $k\geq 3$.
By a time rescaling $\tau=\varepsilon t$,
system (\ref{slow-1}) is changed to
\begin{eqnarray}
\label{fast-1}
\begin{split}
\frac{d x}{d t} &= x'  = f(x,y,\mu,\varepsilon),
\\
\frac{d y}{d t}&= y'  = \varepsilon g(x,y,\mu,\varepsilon).
\end{split}
\end{eqnarray}
For simplification,
let $X_{\varepsilon,\mu}$ denote the vector field of system (\ref{fast-1}).
Clearly, systems (\ref{slow-1}) and (\ref{fast-1}) are equivalent for $\varepsilon\neq 0$.
To obtain the dynamics of system (\ref{slow-1}) or (\ref{fast-1}) for sufficiently small $\varepsilon$,
we  consider the limiting case $\varepsilon=0$.
Then system (\ref{slow-1}) becomes {\it the reduced equation}
\begin{eqnarray}
\label{reduce-1}
\begin{split}
0&=  f(x,y,\mu,0),
\\
\dot y  & = g(x,y,\mu,0),
\end{split}
\end{eqnarray}
and system (\ref{fast-1}) becomes {\it the layer equation}
\begin{eqnarray}
\label{layer-1}
\begin{split}
x'  &= f(x,y,\mu,0),
\\
y'  &= 0,
\end{split}
\end{eqnarray}
For each fixed $\mu$, we observe that
the phase state of system (\ref{reduce-1}) is defined on the set of equilibria  of system (\ref{layer-1}),
that is,
$$\mathcal{C}_{\mu,0}:=\{(x,y)\in \mathbb{R}^{2}: f(x,y,\mu,0)=0\}.$$
This set is called  the {\it critical set}. If it is a submanifold of $\mathbb{R}^{2}$,
then it is called the {\it critical manifold}.
The branches of the set $\mathcal{C}_{\mu,0}$ are called the {\it slow curves}.
By the Fenichel theory \cite{Fenichel-79},
a normally hyperbolic submanifold $\mathcal{M}_{\mu,0}$ (with or without boundary) of
the critical manifold $\mathcal{C}_{\mu,0}$ is perturbed to a {\it slow manifold} $\mathcal{M}_{\mu,\varepsilon}$
near $\mathcal{M}_{\mu,0}$ for sufficiently small $\varepsilon$.
We call $\mathcal{M}_{\mu,0}$ the {\it normally hyperbolic manifold} if
$\partial f/ \partial x\neq 0$ along $\mathcal{M}_{\mu,0}$.
The points in $\mathcal{C}_{\mu,0}$ with $\partial f/ \partial x=0$ are called the {\it contact points},
where the normal hyperbolicity breaks down.
The most generic contact points are {\it jump points},
for which the reduced flow (\ref{reduce-1}) directs towards the contact points.
More degenerate contact points are {\it canard points},
a simple zero of the function $g$ in system (\ref{fast-1}),
which leads to a possibility of periodic orbits in its neighborhood.
Canard points are also called the {\it turning points} in some references.
Geometric analysis of the contact points was initiated in \cite{Dumortieretal-Roussarie-96},
where Dumortier and Roussarie applied  the blow-up technique to study the singularly perturbed van der Pol equation.
Following the pioneering work \cite{Dumortieretal-Roussarie-96} of Dumortier and Roussarie,
many efforts have been devoted to expand the capabilities of this technique.
For example, Krupa and Szmolyan used the technique provided in \cite{Dumortieretal-Roussarie-96}
to extend the slow manifolds of planar singularly perturbed systems near jump points and canard points,
and more results on the blow-up technique and its generalizations are referred to
\cite{Alvarez-11,De-Dumortier-06,De-Dumortier-08,Dumortier-Roussarie-01,Kuehn-15}.

\begin{figure}[!htbp]
\centering
\subfigure[]{
\begin{minipage}[t]{0.24\linewidth}
\centering
\includegraphics[width=1.3in]{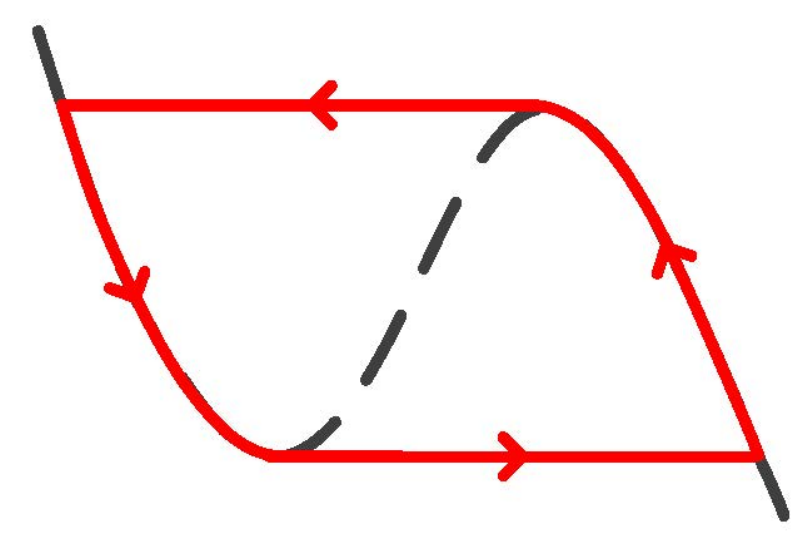}
\end{minipage}
\label{fg-common}
}%
\subfigure[]{
\begin{minipage}[t]{0.24\linewidth}
\centering
\includegraphics[width=1.4in]{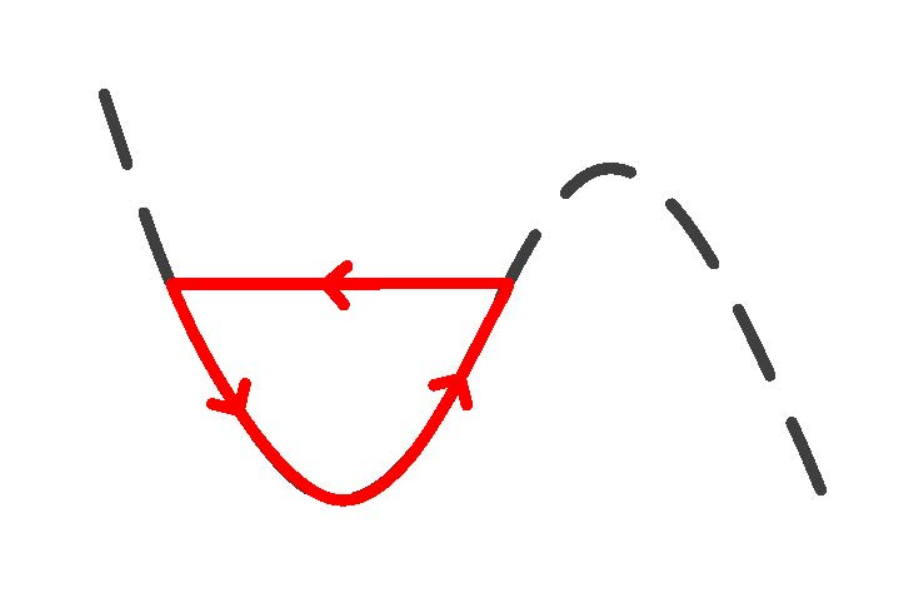}
\end{minipage}
\label{fg-singular-canard-wt}
}%
\subfigure[]{
\begin{minipage}[t]{0.24\linewidth}
\centering
\includegraphics[width=1.4in]{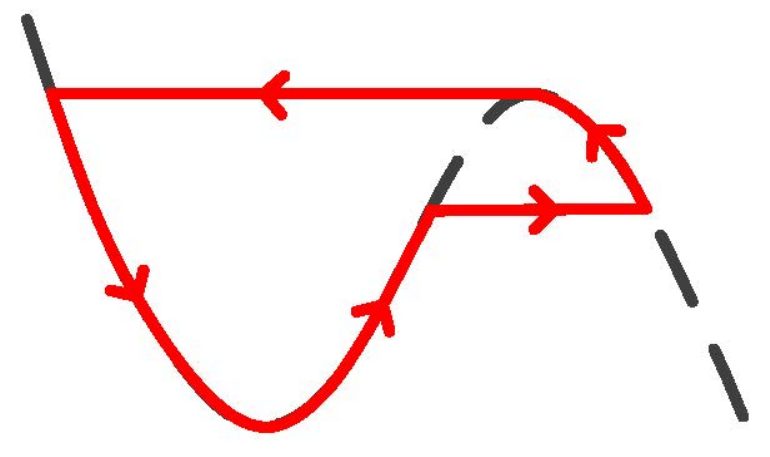}
\end{minipage}
\label{fg-singular-canard-w}
}%
\subfigure[]{
\begin{minipage}[t]{0.24\linewidth}
\centering
\includegraphics[width=1.1in]{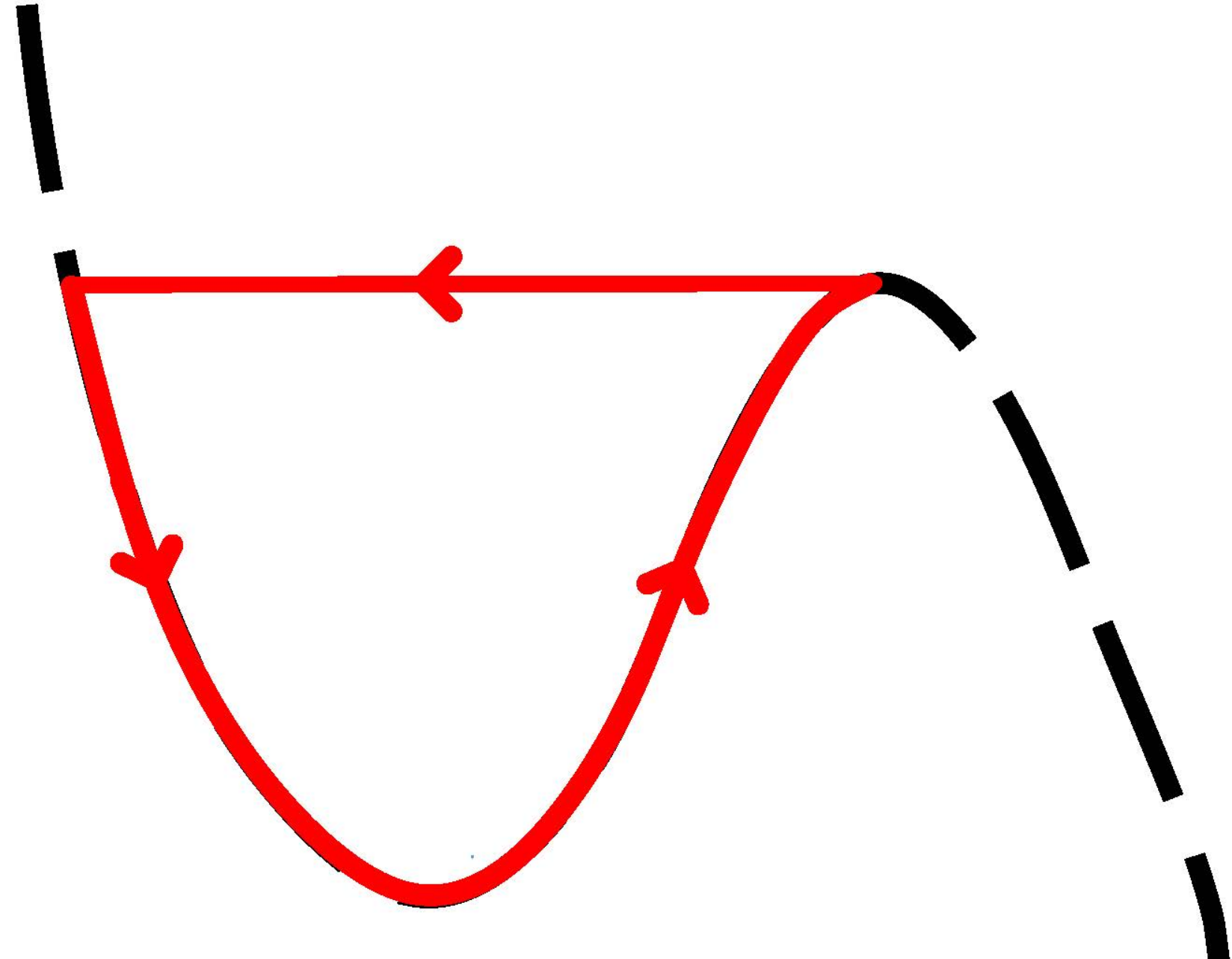}
\end{minipage}
\label{fg-singular-canard-trans}
}%
\centering
\caption{
\ref{fg-common} Common cycle.
 \ref{fg-singular-canard-wt} Canard slow-fast cycle without head.
 \ref{fg-singular-canard-w} Canard slow-fast cycle with head.
 \ref{fg-singular-canard-trans} Transitory canard.
}
\label{fg-slow-fast-cycle}
\end{figure}

Jump points and canard points in planar singularly perturbed  systems with one-parameter layer equations
can lead to relaxation oscillations and canard solutions, respectively.
{\it Relaxtion oscillation} is a periodic orbit which
spends a long time along the slow manifold towards a jump point, jumps from this contact point,
spends a short time parallel to the fast orbits towards another stable slow manifold,
follows the slow manifold again until another jump point is reached,
and finally returns to its starting point via several similarly successive motions \cite{Grasman,Krupa-Szmolyan-01JDE}.
{\it Canards} are the orbits contained in the intersection of
an attracting slow manifold and a repelling slow manifold.
Canards are subject to a {\it generic Hopf breaking mechanism}, that is,
the flow of the layer equation (\ref{layer-1}) has the same direction on a
attracting slow curve and a repelling slow curve which are connected by a generic canard point.
Periodic orbits containing canards are referred to as {\it canard cycles}.
Relaxation oscillations and canard cycles can be both seen as the perturbations of slow-fast cycles
which are closed loops formed by a connected succession of critical manifolds and fast orbits of the layer equations.
More precisely,
relaxation oscillations arise from common cycles
which only contain repelling critical manifolds or attracting critical manifolds (see Figure \ref{fg-common}),
canard cycles from canard slow-fast cycles which contain at least one attracting and one repelling critical manifolds
(see Figures \ref{fg-singular-canard-wt}, \ref{fg-singular-canard-w}  and \ref{fg-singular-canard-trans}).

It is worth mentioning that much work on
relaxation oscillations, canard cycles with head and canard cycles without head bifurcating from  slow-fast cycles
assumed that planar singularly perturbed  systems have precisely one canard points, two jump points or one canard point and one jump point.
However,
there are numerous models of the form (\ref{slow-1}) with
S-shaped critical manifolds possessing two canard points  in real world applications,
such as in a circadian oscillator model based on dimerization and proteolysis of PER and TIM proteins in Drosophila  \cite{chen-duan-li,Tyson-etal-99}, a class of cubic Li\'enard equations with quadratic damping \cite{Dumortier-etal-00}, a predator-prey model of generalized Holling type III \cite{Wang-Zhang-19} and so on.
Stimulated by  it,
we studied the dynamics of system (\ref{fast-1}) with a S-shaped critical manifold,
which has precisely two canard points and one saddle lying on the repelling slow curve.
Under  a certain condition,
there is the simultaneous occurrence of two generic Hopf breaking mechanisms.
\begin{figure}[!htp]
  \centering
  \includegraphics[width=1.7in]{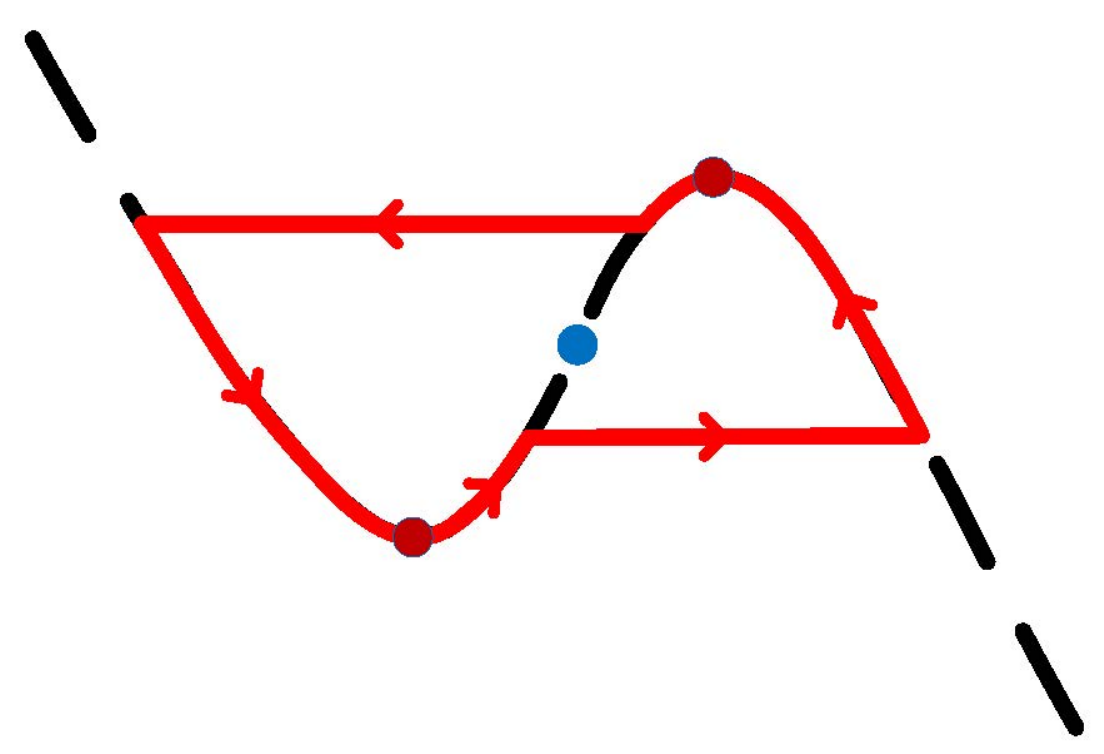}
  \caption{Double canard  slow-fast cycle.}
  \label{fg-double}
\end{figure}
We establish the existence of double canard cycles,
that is, canard cycles include two canards.
This canard cycle can be seen as the limit cycle bifurcating from a double canard slow-fast  cycle
or two-layer canard slow-fast cycles \cite{Dumortier-Roussarie-08},
that is, a slow-fast canard cycle  passes through two layers of fast orbits and
contains two generic Hopf breaking mechanisms at two non-degenerate canard points.
See Figure \ref{fg-double}.
The proofs for these results are based on the normal forms established near the canard points,
the blow-up technique and the extended Melnikov theory obtained by Wechselberger in \cite{Wech-02}.

This paper is organized as follows.
In section \ref{sec-mainresults}
we make some hypotheses and state the main results on the existence of double canard cycles.
In section \ref{sec-normal-form}
we give the normal forms of system (\ref{fast-1}) near canard points.
Section \ref{sec-canard solu} is contributed to investigating the existence of the canards
near canard points by the Melnikov theory.
The proof for the main results is given in section \ref{sec-proof}
and then we apply the main results to cubic Li\'enard equations with quadratic damping in section \ref{sec-lienard}.
In the finial section,
we make some remarks on the limit cycles arising from planar singularly perturbed systems with two canard points.

\section{Main results}
\label{sec-mainresults}
\setcounter{equation}{0}
\setcounter{lemma}{0}
\setcounter{remark}{0}

In this section,
we first introduce some essential hypotheses and state the main results on double canard cycles
arising from system (\ref{slow-1}) with an S-shaped critical manifold possessing two canard points.
We assume that for a fixed $\mu=\mu_{0}=(\lambda_{0},\eta_{1,0},...,\eta_{m,0})$,
the function $f$ satisfies the following hypotheses:
\begin{enumerate}
\item[{\bf (H1)}]
For a fixed $\mu=\mu_{0}$,
there exists a smooth function $\phi$ having precisely two different extreme points $\alpha_{1}$ and $\alpha_{2}$
with $\alpha_{1}<\alpha_{2}$ such that
the critical manifold $\mathcal{C}_{0}:=\mathcal{C}_{\mu_{0},0}$ is represented by
\begin{eqnarray*}
\mathcal{C}_{0}:=\left\{(x,y)\in \mathbb{R}^{2}: y=\phi(x)\right\}.
\end{eqnarray*}

\item[{\bf (H2)}]
At the points $(\alpha_{j},\omega_{j}):=(\alpha_{j},\phi(\alpha_{j}))$, $j=1,2$,
the functions $f$ satisfies the singularities:
\begin{eqnarray}
\label{condit-fold}
f(\alpha_{j},\omega_{j},\mu_{0},0)=0,\ \ \ \ \
\frac{\partial f}{\partial x}(\alpha_{j},\omega_{j},\mu_{0},0)=0,
\end{eqnarray}
and the following non-degenerate conditions:
\begin{eqnarray*}
\frac{\partial^{2} f}{\partial x^{2}}(\alpha_{j},\omega_{j},\mu_{0},0)\neq 0,\ \ \ \
\frac{\partial f}{\partial y}(\alpha_{j},\omega_{j},\mu_{0},0)\neq 0.
\end{eqnarray*}

\item[{\bf (H3)}] Along the critical manifold $\mathcal{C}_{0}$,
the function $f$ satisfies that
\begin{eqnarray*}
&&\frac{\partial f}{\partial x}(x,\phi(x),\mu_{0},0)<0\ \mbox{ for }x\in (-\infty,\alpha_{1})\cup (\alpha_{2},+\infty),\\
&&\frac{\partial f}{\partial x}(x,\phi(x),\mu_{0},0)>0\ \mbox{ for }x\in (\alpha_{1},\alpha_{2}).
\end{eqnarray*}

\end{enumerate}

We remark that by {\bf (H2)} and the {\it Implicit Function Theorem},
there exist an open neighbourhood $U_{\mu}(\mu_{0})\subset\mathbb{R}^{1+m}$ of $\mu=\mu_{0}$
and exactly four $C^{k}$ functions
\begin{eqnarray}
\label{x-expan}
\begin{split}
\tilde{x}_{j}(\mu)
   =&\ \alpha_{j}+\frac{f_{\lambda}f_{xy}-f_{y}f_{x\lambda}}{f_{y}f_{xx}}|_{(\alpha_{j},\omega_{j},\mu_{0},0)}(\lambda-\lambda_{0})\\
    &
   +\sum_{k=1}^{m}\frac{f_{\eta_{k}}f_{xy}-f_{y}f_{x\eta_{k}}}{f_{y}f_{xx}}|_{(\alpha_{j},\omega_{j},\mu_{0},0)}(\eta_{k}-\eta_{k,0})
   +O(\|\mu-\mu_{0}\|^{2}),
\end{split}
\end{eqnarray}
\begin{eqnarray}
\label{y-expan}
\begin{split}
\tilde{y}_{j}(\mu)
 =&\ \omega_{j}+\frac{f_{x}f_{x\lambda}-f_{\lambda}f_{xx}}{f_{y}f_{xx}}|_{(\alpha_{j},\omega_{j},\mu_{0},0)}(\lambda-\lambda_{0})\\
  &\
   +\sum_{k=1}^{m}\frac{f_{x}f_{x\eta_{k}}-f_{\eta_{k}}f_{xx}}{f_{y}f_{xx}}|_{(\alpha_{j},\omega_{j},\mu_{0},0)}(\eta_{k}-\eta_{k,0})
   +O(\|\mu-\mu_{0}\|^{2}),
\end{split}
\end{eqnarray}
such that
\begin{eqnarray}
f(\tilde{x}_{j}(\mu),\tilde{y}_{j}(\mu),\mu,0)=0, \ \ \
\frac{\partial f}{\partial x}(\tilde{x}_{j}(\mu),\tilde{y}_{j}(\mu),\mu,0)=0,\ \ \
\mu\in U_{\mu}(\mu_{0}),\ \ \  j=1,2.
\label{cond-f-1}
\end{eqnarray}
We can choose an appropriate set $U_{\mu}(\mu_{0})$ such that  $(\tilde{x}_{j}(\mu),\tilde{y}_{j}(\mu))$
satisfy the non-degenerate conditions
\begin{eqnarray}
\frac{\partial^{2} f}{\partial x^{2}}(\tilde{x}_{j}(\mu),\tilde{y}_{j}(\mu),\mu,0)\neq 0, \ \ \
\frac{\partial f}{\partial y}(\tilde{x}_{j}(\mu),\tilde{y}_{j}(\mu),\mu,0)\neq 0,
\ \ \ \mu\in U_{\mu}(\mu_{0}),\ \ \  j=1,2.
\label{cond-f-2}
\end{eqnarray}
Then we get two manifolds $\mathcal{L}^{0}_{\mu}$ and $\mathcal{R}^{0}_{\mu}$, which are parameterized by $\mu$ and given by
\begin{eqnarray*}
\mathcal{L}^{0}_{\mu}=\left\{(\tilde{x}_{1}(\mu),\tilde{y}_{1}(\mu))\in\mathbb{R}^{2}: \mu\in U_{\mu}(\mu_{0})\right\},\ \ \
\mathcal{R}^{0}_{\mu}=\left\{(\tilde{x}_{2}(\mu),\tilde{y}_{2}(\mu))\in\mathbb{R}^{2}: \mu\in U_{\mu}(\mu_{0})\right\}.
\end{eqnarray*}
We refer to $\mathcal{L}^{0}_{\mu}$ and $\mathcal{R}^{0}_{\mu}$ as the {\it contact point manifolds}.
Without loss of generality, assume that $U_{\mu}(\mu_{0})=\mathbb{R}^{1+m}$.

Under the above hypotheses,
we see that $(\alpha_{j},\omega_{j})$ are both contact points.
Our interest is to study the case that they are both canard points.
Then we further make the following hypotheses on the function $g$.
\begin{enumerate}
\item[{\bf (H4)}]
For $\mu=\mu_{0}$, the function $g$ satisfies that
\begin{eqnarray*}
g(\alpha_{j},\omega_{j},\mu_{0},0)=0,\ \ \
\frac{\partial g}{\partial x}(\alpha_{j},\omega_{j},\mu_{0},0)\neq 0, \ \ \ j=1,2,
\end{eqnarray*}
and in the extended space $\{(x,y,\lambda)\in\mathbb{R}^{3}\}$,
the curves $\mathcal{L}^{0}_{\mu}$ and $\mathcal{R}^{0}_{\mu}$  transversely intersect the manifold given by $g(x,y,\lambda,\eta_{0},0)=0$
at $(\alpha_{j},\omega_{j},0)$,
that is,
\begin{eqnarray*}
G_{j}:=
\frac{\partial g}{\partial x}(\alpha_{j},\omega_{j},\lambda_{0},\eta_{0},0)\cdot \frac{\partial\tilde{x}_{j}}{\partial \lambda}(\mu_{0})
+\frac{\partial g}{\partial y}(\alpha_{j},\omega_{j},\lambda_{0},\eta_{0},0)\cdot \frac{\partial\tilde{y}_{j}}{\partial \lambda}(\mu_{0})
+\frac{\partial g}{\partial \lambda}(\alpha_{j},\omega_{j},\lambda_{0},\eta_{0},0)
\neq 0.
\end{eqnarray*}

\item[{\bf (H5)}]
For $\mu=\mu_{0}$,
system (\ref{slow-1}) has precisely one equilibrium $E_{0}:=(x_{m},y_{m})$ on the section $M$ which is a saddle,
the function $\phi$ respectively reaches its minimum  and maximum  values at $x=\alpha_{1}$ and $x=\alpha_{2}$,
and the slow motions governed by
\begin{eqnarray*}
\phi'(x)\dot x=\phi'(x)\frac{dx}{d\tau}=g(x,\phi(x),\mu_{0},0),
\end{eqnarray*}
satisfy $\dot x>0$ for $x<\alpha_{1}$ and $\dot x<0$ for $x>\alpha_{2}$.
\end{enumerate}

\begin{figure}[!htbp]
\centering
\subfigure[]{
\begin{minipage}[t]{0.45\linewidth}
\centering
\includegraphics[width=5.8cm]{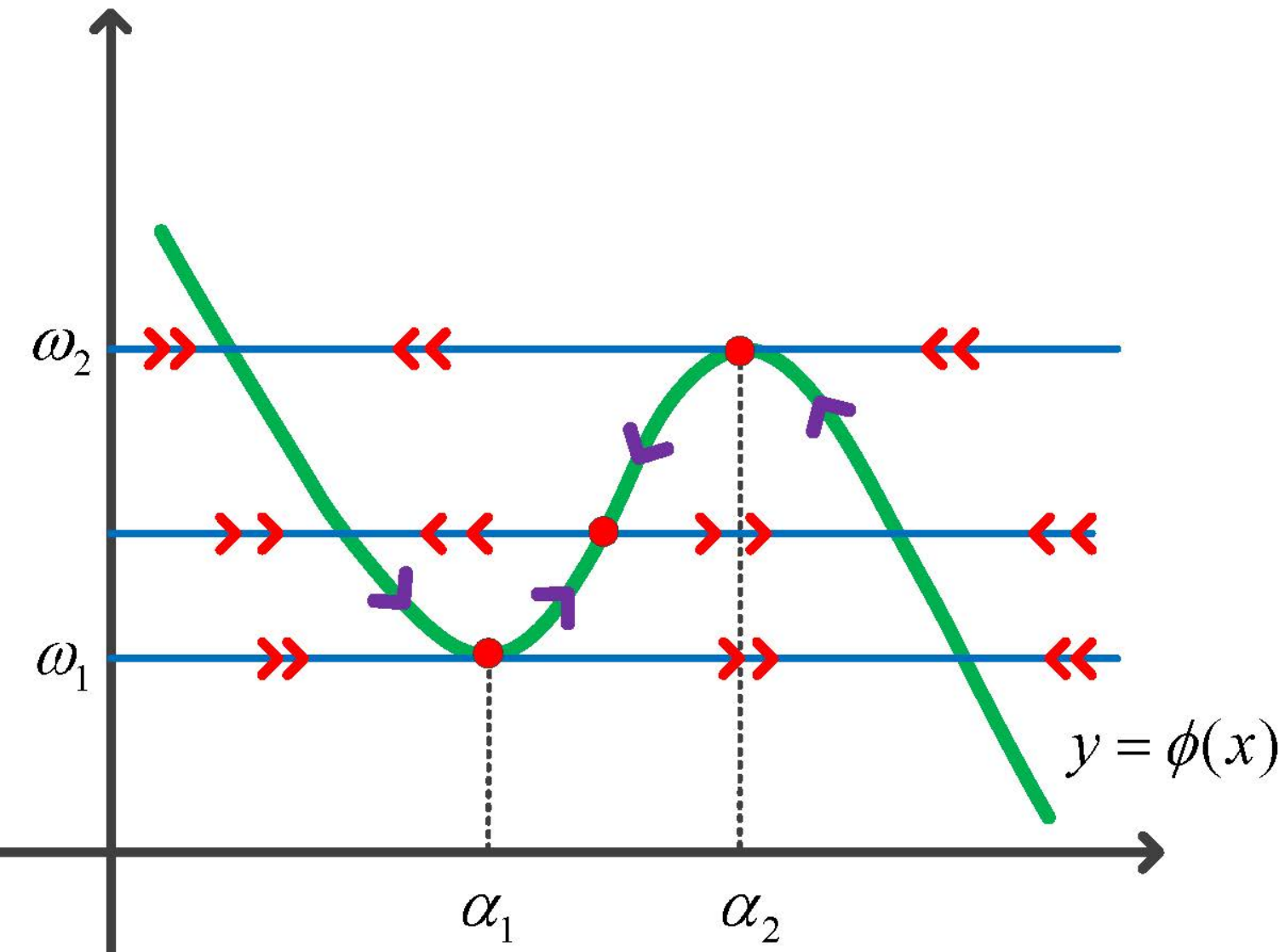}
\end{minipage}
\label{fg-S-F-2}
}%
\subfigure[]{
\begin{minipage}[t]{0.45\linewidth}
\centering
\includegraphics[width=5.8cm]{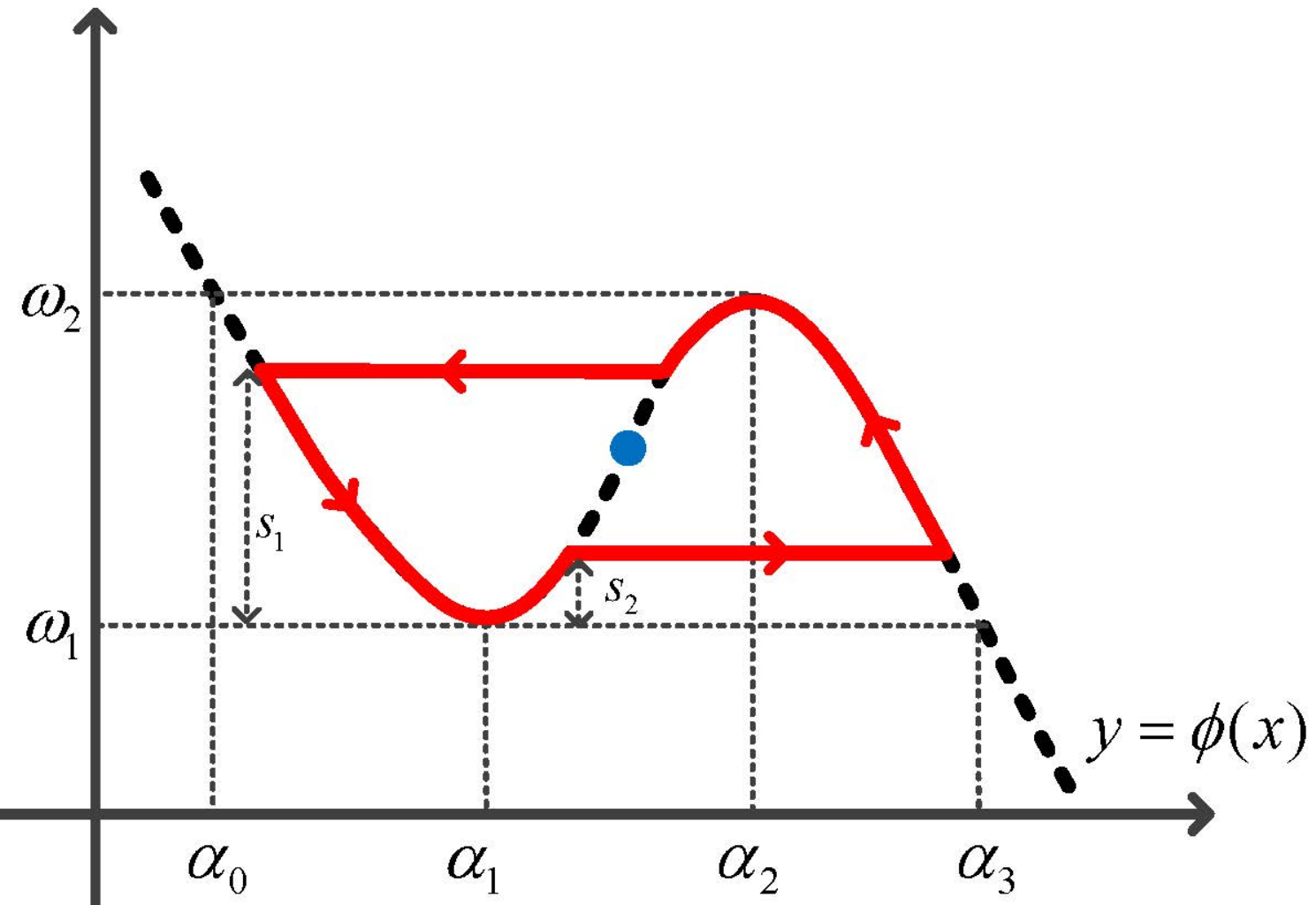}
\end{minipage}
\label{fg-S-Fig-two-canards}
}%
\centering
\caption{\ref{fg-S-F-2} The behaviors of the limiting systems (\ref{reduce-1}) and (\ref{layer-1}).
\ref{fg-S-Fig-two-canards} Double canard slow-fast cycle.}
\end{figure}

Under the hypotheses {\bf (H1)}-{\bf (H5)},
for $\mu=\mu_{0}$ both contact points are canard points.
The dynamics of the limiting systems (\ref{reduce-1}) and (\ref{layer-1}) are shown in Figure \ref{fg-S-F-2}.
Our goal is to study the double canard cycles
arising from system (\ref{slow-1}) with {\bf (H1)}-{\bf (H5)}.
To establish the existence of double canard cycles,
we construct  double canard slow-fast cycles in the following way.
For each $s\in (0,\omega_{2}-\omega_{1})$,
let the constants $\alpha_{L}(s)$, $\alpha_{M}(s)$ and $\alpha_{R}(s)$
satisfy $\alpha_{L}(s)<\alpha_{1}<\alpha_{M}(s)<\alpha_{2}<\alpha_{R}(s)$
and $\phi(\alpha_{L}(s))=\phi(\alpha_{M}(s))=\phi(\alpha_{R}(s))=\phi(\omega_{1}+s)$.
For a pair of real constants $(s_{1},s_{2})$ in  $\Omega:=\left\{0<s_{2}<y_{m}-\omega_{1}<s_{1}<\omega_{2}-\omega_{1}\right\}$,
define a slow-fast cycle $\Gamma(s_1,s_{2})$ by
(see Figure \ref{fg-S-Fig-two-canards})
\begin{eqnarray}
\label{df-butterfly}
\begin{split}
\Gamma(s_1,s_{2})=& \left\{(x,y)\in\mathbb{R}^{2}: y=\phi(x),\ \alpha_{L}(s_{1})\leq x\leq \alpha_{M}(s_{2})\right\}\\
          &  \cup \left\{(x,\omega_{1}+s_{2})\in\mathbb{R}^{2}:\alpha_{M}(s_{2})\leq x\leq \alpha_{R}(s_{2})\right\}\\
          &  \cup \left\{(x,y)\in\mathbb{R}^{2}: y=\phi(x),\ \alpha_{M}(s_{1})\leq x\leq \alpha_{R}(s_{2})\right\}\\
          &  \cup \left\{(x,\omega_{1}+s_{1})\in\mathbb{R}^{2}:\alpha_{L}(s_{1})\leq x\leq \alpha_{M}(s_{1})\right\}.
\end{split}
\end{eqnarray}
The limit cycles bifurcating from $\Gamma(s_1,s_{2})$ is closely related to slow divergence integrals.
Here
for any $x_{1}$ and $x_{2}$ in $\mathbb{R}$ the slow divergence integral $I(x_{1},x_{2})$ along the critical manifold is in the form
(see, for instance, \cite{De-etal-11,De-etal-14,Dumortier-11,Krupa-Szmolyan-01JDE,Li-Zhu-13,Li-Lu-14})
\begin{eqnarray*}
I(x_{1},x_{2})=\int_{x_{1}}^{x_{2}}
   \frac{\partial f}{\partial x}(x,\phi(x),\mu_{0},0)\cdot\frac{\phi'(x)}{g(x,\phi(x),\mu_{0},0)}\,dx.
\end{eqnarray*}
Four slow divergence integrals associated with $\Gamma(s_1,s_{2})$ are given by
\begin{eqnarray*}
\mathcal{I}_{1}(s_{1})=I(\alpha_{L}(s_{1}),\alpha_{1}), \ \
\mathcal{I}_{2}(s_{2})=I(\alpha_{1},\alpha_{M}(s_{2})), \ \
\mathcal{I}_{3}(s_{2})=I(\alpha_{R}(s_{2}),\alpha_{2}), \ \
\mathcal{I}_{3}(s_{2})=I(\alpha_{2},\alpha_{M}(s_{1})).
\end{eqnarray*}
Before stating the main results,
we define some important constants which are useful in the subsequent proof.
Let the constants $d_{j,r_{2}}$, $d_{j,\lambda_{2}}$ and $d_{j,\eta_{i}}$ be in the form
\begin{eqnarray}
\label{df-coeff}
d_{j,r_{2}}=-\frac{\sqrt{2\pi}}{8}(4a_{j,1}-a_{j,2}+3a_{j,3}-2a_{j,4}+2a_{j,5}),\ \ \
d_{j,\lambda_{2}}=-\sqrt{2\pi},  \ \ \ d_{j,\eta_{2,i}}=\sqrt{2\pi}a_{j,5+i}, i=1,...,m.
\end{eqnarray}
where  $a_{j,i}$ are given by
\begin{eqnarray*}
&&a_{j,1}=\frac{\partial \phi_{j,3}}{\partial x}(0), \ \ \
a_{j,2}=\frac{\partial \phi_{j,1}}{\partial x}(0),\ \ \
a_{j,3}=\frac{\partial \phi_{j,2}}{\partial x}(0),\\
&&a_{j,4}=\frac{\partial \phi_{j,4}}{\partial x}(0),\ \ \
a_{j,5}=\phi_{j,6}(0),\ \ \
a_{j,5+i}=\phi_{j,6+i}(0), \ \ \ j=1,2, \ \ \ i=1,...,m,
\end{eqnarray*}
and $\phi_{j,i}=\phi_{i}$ whose lengthy expressions are given as in Lemma \ref{lm-normal-form} for $(\alpha,\omega)=(\alpha_{j},\omega_{j})$.
The main theorem are stated in the following.

\renewcommand\thetheorem{A}
\begin{theorem}
Suppose that system (1.1) satisfies {\bf (H1)}-{\bf (H5)}.
Let the functions $\beta_{j}$ be given by
$\beta_{j}(\mu)=g(\tilde{x}_{j}(\mu),\tilde{y}_{j}(\mu),\mu,0)$ for $\mu\in \mathbb{R}^{1+m}$ and $j=1,2$.
Assume that the curves
$\left\{(s_{1},s_{2})\in\Omega: \mathcal{I}_{1}(s_{1}) +\mathcal{I}_{2}(s_{2})=0 \right\}$
and $\left\{(s_{1},s_{2})\in\Omega: \mathcal{I}_{3}(s_{2}) +\mathcal{I}_{4}(s_{1})=0 \right\}$
transversally intersects at the point $(s_{1}^{*},s_{2}^{*})\in\Omega$,
and for a certain $i\in \{1,...,m\}$,
\begin{eqnarray}
\label{con-rk-2}
{\rm Rank}
\left(
\begin{array}{ll}
-f_{xx}(P_{1})(2g_{x}(P_{1}))^{-1}G_{1}d_{1,\lambda_{2}} & d_{1,\eta_{2,i}} \\
-f_{xx}(P_{2})(2g_{x}(P_{2}))^{-1}G_{2}d_{2,\lambda_{2}} & d_{2,\eta_{2,i}}
\end{array}
\right)
={\rm Rank}
\left(
\begin{array}{ll}
\frac{\partial \beta_{1}}{\partial \lambda}(\mu_{0}) & \frac{\partial \beta_{1}}{\partial \eta_{i}}(\mu_{0})  \\
\frac{\partial \beta_{2}}{\partial \lambda}(\mu_{0})  & \frac{\partial \beta_{2}}{\partial \eta_{i}}(\mu_{0})
\end{array}
\right)
=2.
\end{eqnarray}
Then there exists a sufficiently small $\varepsilon_{0}$ and  two  continuous functions $\lambda(\varepsilon)$ and $\eta_{i}(\varepsilon)$
of the form
\begin{eqnarray}
\label{df-lambda-eta}
\begin{split}
\lambda(\varepsilon)&=&\!\!\!\lambda_{0}+(A_{1}B_{2}-A_{2}B_{1})^{-1}(B_{2}C_{1}-B_{1}C_{2})\varepsilon+O(\varepsilon^{3/2}),\\
\eta_{i}(\varepsilon)&=&\!\!\!\eta_{i,0}+(A_{1}B_{2}-A_{2}B_{1})^{-1}(A_{1}C_{2}-A_{2}C_{1})\varepsilon+O(\varepsilon^{3/2}),
\end{split}
\end{eqnarray}
for $0\leq \varepsilon<\varepsilon_{0}$,
where $A_{j}$, $B_{j}$ and $C_{j}$ are given by
\begin{eqnarray*}
A_{j}=-\frac{f_{xx}(P_{j})G_{j}}{2g_{x}(P_{j})}d_{j,r_{2}},\ \ \
B_{j}=d_{j,\eta_{2,i}}, \ \ \
C_{j}=f_{y}(P_{j})g_{x}(P_{j})d_{j,r_{2}}+\frac{f_{xx}(P_{j})g_{\varepsilon}(P_{j})}{2g_{x}(P_{j})}d_{j,\lambda_{2}},
\ \  \ j=1,2,
\end{eqnarray*}
such that system (1.1) with $(\varepsilon,\mu)=(\varepsilon,\lambda(\varepsilon),\eta_{1,0},...,\eta_{i-1,0},\eta_{i}(\varepsilon),\eta_{i+1,0},...,\eta_{m,0})$
has a limit cycle $\Gamma_{\varepsilon}(s_1,s_{2})$ in a small neighborhood of $\Gamma(s_1,s_{2})$,
and  $\Gamma_{\varepsilon}(s_1,s_{2})\to \Gamma(s_1,s_{2})$ as $\varepsilon\to 0$ in the sense of Hausdorff distance.
Furthermore, there exist three hyperbolic limit cycles arising from $\Gamma_{\varepsilon}(s_1,s_{2})$ under some suitable conditions.
\label{thm-exist}
\end{theorem}
%

\section{Normal forms near canard points}
\label{sec-normal-form}

In this section, we consider the normal forms of system (\ref{fast-1}) near the canard points.
Let the functions $\tilde{g}_{j}$ associated with $(\alpha_{j},\omega_{j})$ be  defined by
\begin{eqnarray*}
\tilde{g}_{j}(x,y,\mu,\varepsilon)=g\left(x+\tilde{x}_{j}(\mu+\mu_{0}), y+\tilde{y}_{j}(\mu+\mu_{0}), \mu+\mu_{0}, \varepsilon\right),
\ \ \ j=1,2.
\end{eqnarray*}
By {\bf (H4)}
we obtain that the functions $\tilde{g}_{j}$ satisfy the following:
\begin{eqnarray}
\tilde{g}_{j}(0)=0,\ \ \
\frac{\partial \tilde{g}_{j}}{\partial x}(0)\neq 0, \ \ \
\frac{\partial \tilde{g}_{j}}{\partial \lambda}(0)\neq 0, \ \ \ j=1,2.
\label{condt-g}
\end{eqnarray}
Whenever there is no confusion, the zero vector is denoted by $0$.
Then by the {\it Implicit Function Theorem},
there exists an open neighbourhood $U_{\varepsilon}(0)\subset\mathbb{R}$ of $\varepsilon=0$
and exactly two $C^{k}$ functions
\begin{eqnarray}
\label{lambda-expan}
\tilde{\lambda}_{j}(\varepsilon)
   =-\frac{\partial\tilde{g}_{j}}{\partial \varepsilon}(0)\cdot
   \left(\frac{\partial\tilde{g}_{j}}{\partial \lambda}(0)\right)^{-1}\varepsilon+O(\varepsilon^{2}),
\end{eqnarray}
such that
\begin{eqnarray}
\tilde{g}_{j}(0,0,\tilde{\lambda}_{j}(\varepsilon),0,\varepsilon)=0,
\ \ \ \ \varepsilon\in U_{\varepsilon}(0), \ \ \ j=1,2.
\label{cond-g-1}
\end{eqnarray}
Without loss of generality we assume that $U_{\varepsilon}(0)=\mathbb{R}$.
The normal forms of system (\ref{fast-1}) near the canard points are given in the next lemma.

\begin{lemma}
\label{lm-normal-form}
Suppose that the functions $f$ and $g$ in system (\ref{fast-1}) respectively
satisfy {\bf (H2)} and {\bf (H4)} at $(\alpha_{j},\omega_{j})$, $j=1,2$.
For each $j=1,2$,
let the functions $\bar{f}$ and $\bar{g}$ be in the form
\begin{eqnarray*}
\bar{f}(x,y,\lambda,\eta,\varepsilon)
 \!\!\!&=&\!\!\!
 f(x+\tilde{x}((\lambda+\tilde{\lambda}(\varepsilon),\eta)+\mu_{0}),y+\tilde{y}((\lambda+\tilde{\lambda}(\varepsilon),\eta)+\mu_{0}),
   (\lambda+\tilde{\lambda}(\varepsilon),\eta)+\mu_{0},\varepsilon),\\
\bar{g}(x,y,\lambda,\eta,\varepsilon)
 \!\!\!&=&\!\!\!
 g(x+\tilde{x}((\lambda+\tilde{\lambda}(\varepsilon),\eta)+\mu_{0}),y+\tilde{y}((\lambda+\tilde{\lambda}(\varepsilon),\eta)+\mu_{0}),
   (\lambda+\tilde{\lambda}(\varepsilon),\eta)+\mu_{0},\varepsilon),
\end{eqnarray*}
where the $C^{k}$ functions $\tilde{x}=\tilde{x}_{j}$, $\tilde{y}=\tilde{y}_{j}$
and $\tilde{\lambda}=\tilde{\lambda}_{j}$ respectively have the expansions (\ref{x-expan}), (\ref{y-expan}) and (\ref{lambda-expan}).
Then near the point $(\alpha,\omega)=(\alpha_{j},\omega_{j})$,
system (\ref{fast-1}) can be changed into the form
\begin{eqnarray}
\label{fast-normal-1}
\begin{split}
 x'  &=
    -y\left(1+\phi_{1}(x,y,\lambda,\eta,\varepsilon)\right)
    +x^{2}\left(1+\phi_{2}(x,y,\lambda,\eta,\varepsilon)\right)+\varepsilon\phi_{3}(x,y,\lambda,\eta,\varepsilon),
\\
y'  & =
    \varepsilon \left(\zeta x\left(1+\phi_{4}(x,y,\lambda,\eta,\varepsilon)\right)
    -\lambda\left(1+\phi_{5}(x,y,\lambda,\eta,\varepsilon)\right)+y\phi_{6}(x,y,\lambda,\varepsilon)
    +\sum_{i=1}^{m}\eta_{i}\phi_{6+i}(x,y,\lambda,\eta,\varepsilon)
    \right),
\end{split}
\end{eqnarray}
where $\zeta=\pm 1$, the $C^{k}$  functions $\phi_{i}$ are defined by
\begin{eqnarray*}
\phi_{1}(x,y,\mu,\varepsilon)\!\!\!&=&\!\!\!
             \left(\bar{f}_{y}(0)\right)^{-1}\hat{\phi}_{1}\circ\mathcal{T}(x,y,\mu,\varepsilon),\\
\phi_{2}(x,y,\mu,\varepsilon)\!\!\!&=&\!\!\!
             2\left(\bar{f}_{xx}(0)\right)^{-1}\hat{\phi}_{2}\circ\mathcal{T}(x,y,\mu,\varepsilon)),\\
\phi_{3}(x,y,\mu,\varepsilon)\!\!\!&=&\!\!\!
           -\zeta\bar{f}_{xx}(0)\left(2\bar{f}_{y}(0)\bar{g}_{x}(0)\right)^{-1}
           \left(\bar{f}_{\varepsilon}(0)+\hat{\phi}_{3}\circ\mathcal{T}(x,y,\mu,\varepsilon)\right),\\
\phi_{4}(x,y,\mu,\varepsilon)\!\!\!&=&\!\!\!
     \left(\bar{g}_{x}(0)\right)^{-1}\hat{\phi}_{4}\circ\mathcal{T}(x,y,\mu,\varepsilon),\\
\phi_{5}(x,y,\mu,\varepsilon)\!\!\!&=&\!\!\!
     \left(\bar{g}_{\lambda}(0)\right)^{-1}\hat{\phi}_{5}\circ\mathcal{T}(x,y,\mu,\varepsilon),\\
\phi_{6}(x,y,\mu,\varepsilon)\!\!\!&=&\!\!\!
    -\zeta \left(\bar{f}_{y}(0)\bar{g}_{x}(0)\right)^{-1}\left(\bar{g}_{y}(0)+\hat{\phi}_{6}\circ\mathcal{T}(x,y,\mu,\varepsilon)\right),\\
\phi_{6+j}(x,y,\mu,\varepsilon) \!\!\!&=&\!\!\!
    \zeta  \bar{f}_{xx}(0) (2\bar{g}_{x}(0))^{-1}\left(\bar{g}_{\eta_{6+j}}(0)+\hat{\phi}_{6+j}\circ\mathcal{T}(x,y,\mu,\varepsilon)\right),
    \ \ \ j=1,...,m,
\end{eqnarray*}
the function $\hat{\phi}_{j}$  have the form
\begin{eqnarray*}
\hat{\phi}_{1}(x,y,\mu,\varepsilon)
       \!\!\!&=&\!\!\!
           -\bar{f}_{y}(0)+\int_{0}^{1}D_{2}\bar{f}(u x, u y, \mu, u \varepsilon) du
           +x\int_{0}^{1}\int_{0}^{1}u D_{12}\bar{f}(uv x, uv y, \mu,uv\varepsilon)dv  d u,\\
\hat{\phi}_{2}(x,y,\mu,\varepsilon)
        \!\!\!&=&\!\!\!
         -\frac{1}{2}\bar{f}_{xx}(0)+\int_{0}^{1}\int_{0}^{1}u D_{11}\bar{f}(uv x, uv y, \mu,uv\varepsilon)dv  d u,\\
\hat{\phi}_{3}(x,y,\mu,\varepsilon)
        \!\!\!&=&\!\!\!
           -\bar{f}_{\varepsilon}(0)+\int_{0}^{1}D_{4+m}\bar{f}(u x, u y, \mu, u \varepsilon) du
           +x\int_{0}^{1}\int_{0}^{1}u D_{14+m}\bar{f}(uv x, uv y, \mu,uv\varepsilon)dv  d u,\\
\hat{\phi}_{4}(x,y,\mu,\varepsilon)
       \!\!\!&=&\!\!\!
           -\bar{g}_{x}(0)+\int_{0}^{1}D_{1}\bar{g}(u x, u y, u \lambda,u \eta, \varepsilon) du,\\
\hat{\phi}_{5}(x,y,\mu,\varepsilon)
        \!\!\!&=&\!\!\!
         -\bar{g}_{\lambda}(0)+\int_{0}^{1}D_{3}\bar{g}(u x, u y, u \lambda,u \eta,  \varepsilon) du,\\
\hat{\phi}_{6}(x,y,\mu,\varepsilon)
        \!\!\!&=&\!\!\!
          -\bar{g}_{y}(0)+\int_{0}^{1}D_{2}\bar{g}(u x, u y, u \lambda,u \eta,  \varepsilon) du,\\
\hat{\phi}_{6+j}(x,y,\mu,\varepsilon)
        \!\!\!&=&\!\!\!
          -\bar{g}_{\eta_{j}}(0)+\int_{0}^{1}D_{3+k}\bar{g}(u x, u y, u \lambda, u \eta,  \varepsilon) du,
          \ \ \ j=1,...,m,
\end{eqnarray*}
the transformation $\tilde{\mathcal{T}}:\mathbb{R}^{4+m}\to\mathbb{R}^{4+m}$ is in the form
\begin{eqnarray*}
\tilde{\mathcal{T}}(x,y,\lambda,\eta,\varepsilon)
      =\left(\frac{2}{\bar{f}_{xx}(0)}x,\,
              -\frac{2}{\bar{f}_{xx}(0)\bar{f}_{y}(0)}y,\,
              -\frac{2\zeta \bar{g}_{x}(0)}{\bar{f}_{xx}(0)\bar{g}_{\lambda}(0)}\lambda,\,
              \eta,\,
              -\frac{\zeta }{\bar{f}_{y}(0)\bar{g}_{x}(0)}\varepsilon
      \right),
\end{eqnarray*}
and $D_{j}$ denotes the {\it j-}th partial derivative with respect to the {\it j-}th variable,
and $D_{ij}$ with $i=1$, $j=1,2,4+m$, are in the form $D_{ij}=D_{j}\circ D_{i}$.
\end{lemma}
{\bf Proof.}
It suffices to study the following system
\begin{eqnarray}
\label{normal-1}
\begin{split}
\frac{d x}{d t} &= x'  = \bar{f}(x,y,\mu,\varepsilon),
\\
\frac{d y}{d t}&= y'  = \varepsilon \bar{g}(x,y,\mu,\varepsilon).
\end{split}
\end{eqnarray}
From (\ref{cond-f-1}), (\ref{cond-f-2}), (\ref{condt-g}) and (\ref{cond-g-1}),
the following hold:
\begin{eqnarray}
\!\!\!\!\!\!\!\!&&\bar{f}(0,0,\mu,0)=0, \ \ \
\frac{\partial \bar{f}}{\partial x}(0,0,\mu,0)=0,\ \ \
\frac{\partial^{2} \bar{f}}{\partial x^{2}}(0,0,\mu,0)\neq 0, \ \ \
\frac{\partial \bar{f}}{\partial y}(0,0,\mu,0)\neq 0, \ \ \ \forall\,\mu\in \mathbb{R}^{1+m},\\
\label{cond-f}
\!\!\!\!\!\!\!\!&&\bar{g}(0)=0,\ \ \
\frac{\partial \bar{g}}{\partial x}(0)\neq 0, \ \ \
\frac{\partial \bar{g}}{\partial \lambda}(0)\neq 0, \ \ \
\bar{g}(0,0,0,0,\varepsilon)=0, \ \ \ \forall\,\varepsilon\in\mathbb{R}.
\label{cond-g}
\end{eqnarray}
Since the function $\bar{f}$ satisfies $\bar{f}(0,0,\mu,0)=0$ for each $\mu\in \mathbb{R}^{1+m}$,
then by \cite[Lemma 2.1, p.5]{Milnor-63} we obtain
\begin{eqnarray*}
\bar{f}(x,y,\mu,\varepsilon)
       \!\!\!&=&\!\!\!
          \int_{0}^{1}\frac{\partial \bar{f}}{\partial u}(u x, u y, \mu, u \varepsilon) du   \\
       \!\!\!&=&\!\!\!
          x\int_{0}^{1}D_{1}\bar{f}(u x, u y, \mu, u \varepsilon) du
          +y\int_{0}^{1}D_{2}\bar{f}(u x, u y, \mu, u \varepsilon) du\\
           \!\!\!& &\!\!\!
          +\varepsilon \int_{0}^{1}D_{4+m}\bar{f}(u x, u y, \mu, u \varepsilon) du.
\end{eqnarray*}
Similarly, we have
\begin{eqnarray*}
D_{1}\bar{f}(u x, u y, \mu, u \varepsilon)
        \!\!\!&=&\!\!\!
   D_{1}\bar{f}(u x, u y, \mu, u \varepsilon)-D_{1}\bar{f}(0, 0, \mu, 0)\\
       \!\!\!&=&\!\!\!
          u x\int_{0}^{1}D_{11}\bar{f}(uv x, uv y, \mu,uv\varepsilon) dv
          +u y\int_{0}^{1}D_{12}\bar{f}(uv x, uv y, \mu,uv\varepsilon) dv\\
           \!\!\!& &\!\!\!
          +u\varepsilon \int_{0}^{1}D_{14+m}\bar{f}(uv x, uv y, \lambda,uv\varepsilon) dv.
\end{eqnarray*}
Thus, the function $\bar{f}$ can be written as the form
\begin{eqnarray}
\label{exp-f-1}
\begin{split}
\bar{f}(x,y,\mu,\varepsilon)\,=&\
     y\left(\bar{f}_{y}(0)+\hat{\phi}_{1}(x,y,\mu,\varepsilon)\right)+
       x^{2}\left(\frac{1}{2}\bar{f}_{xx}(0)+\hat{\phi}_{2}(x,y,\mu,\varepsilon)\right)
                           \\
       & +\varepsilon\left(\bar{f}_{\varepsilon}(0)+\hat{\phi}_{3}(x,y,\mu,\varepsilon)\right),
\end{split}
\end{eqnarray}
where $\hat{\phi}_{i}$ are defined as in this lemma.

Similarly, note that $\bar{g}(0,0,0,0,\varepsilon)=0$ for $\varepsilon\in U_{\varepsilon}(0)$,
then $\bar{g}$ can be written as the form
\begin{eqnarray}
\label{exp-g-1}
\begin{split}
\bar{g}(x,y,\lambda,\eta,\varepsilon)=& \ \bar{g}(x,y,\lambda,\eta,\varepsilon)-\bar{g}(0,0,0,0,\varepsilon)\\
=&\   x\left(\bar{g}_{x}(0)+\hat{\phi}_{4}(x,y,\lambda,\eta,\varepsilon)\right)
      +\lambda\left(\bar{g}_{\lambda}(0)+\hat{\phi}_{5}(x,y,\lambda,\eta,\varepsilon)\right)\\
      &\ +y\left(\bar{g}_{y}(0)+\hat{\phi}_{6}(x,y,\lambda,\eta,\varepsilon)\right)
      +\sum_{j=1}^{m}\eta_{j}\left(\bar{g}_{\eta_{j}}(0)+\hat{\phi}_{6+j}(x,y,\lambda,\eta,\varepsilon)\right),
\end{split}
\end{eqnarray}
where the functions $\hat{\phi}_{i}$, $i=4,...,6+m$, are defined as in this lemma.
Clearly, the functions $\hat{\phi}_{i}$ and $\hat{\psi}_{j}$, $i=1,2,3$, $j=1,...,3+m$,
are $C^{k}$ and satisfy  $\hat{\phi}_{i}(0)=0$ and $\hat{\psi}_{j}(0)=0$.
Then by taking the transformation $\tilde{\mathcal{T}}$,
we obtain the normal form (\ref{fast-normal-1}).
Therefore, the proof is now complete.
\hfill$\Box$
%

\section{Existence of canards}
\label{sec-canard solu}
\setcounter{equation}{0}
\setcounter{lemma}{0}
\setcounter{theorem}{0}
\setcounter{remark}{0}

\begin{figure}[!htbp]
\centering
\subfigure[]{
\begin{minipage}[t]{0.45\linewidth}
\centering
\includegraphics[width=3.6cm]{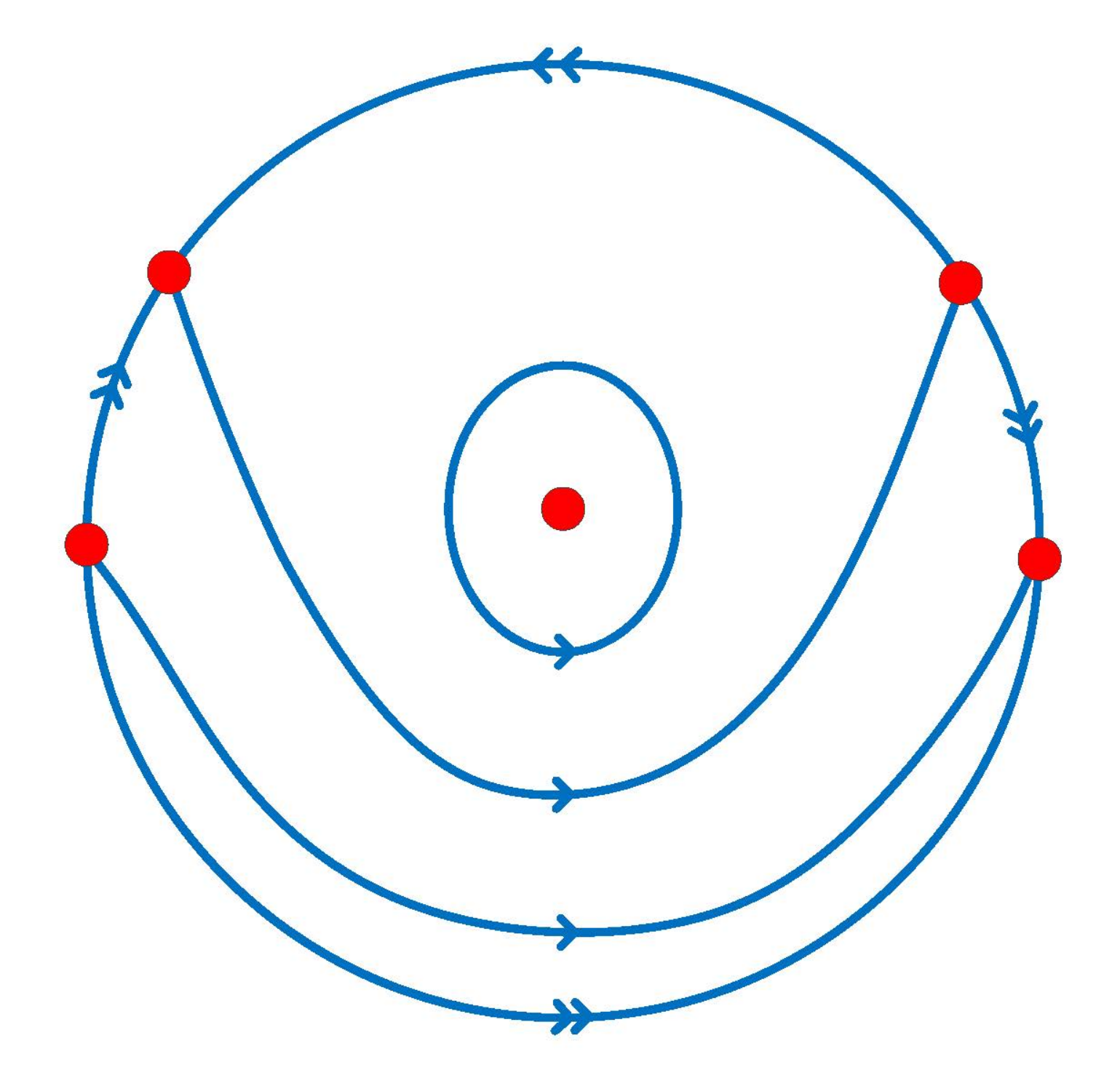}
\end{minipage}
\label{fg-poincare-sphere}
}%
\subfigure[]{
\begin{minipage}[t]{0.45\linewidth}
\centering
\includegraphics[width=3.4cm]{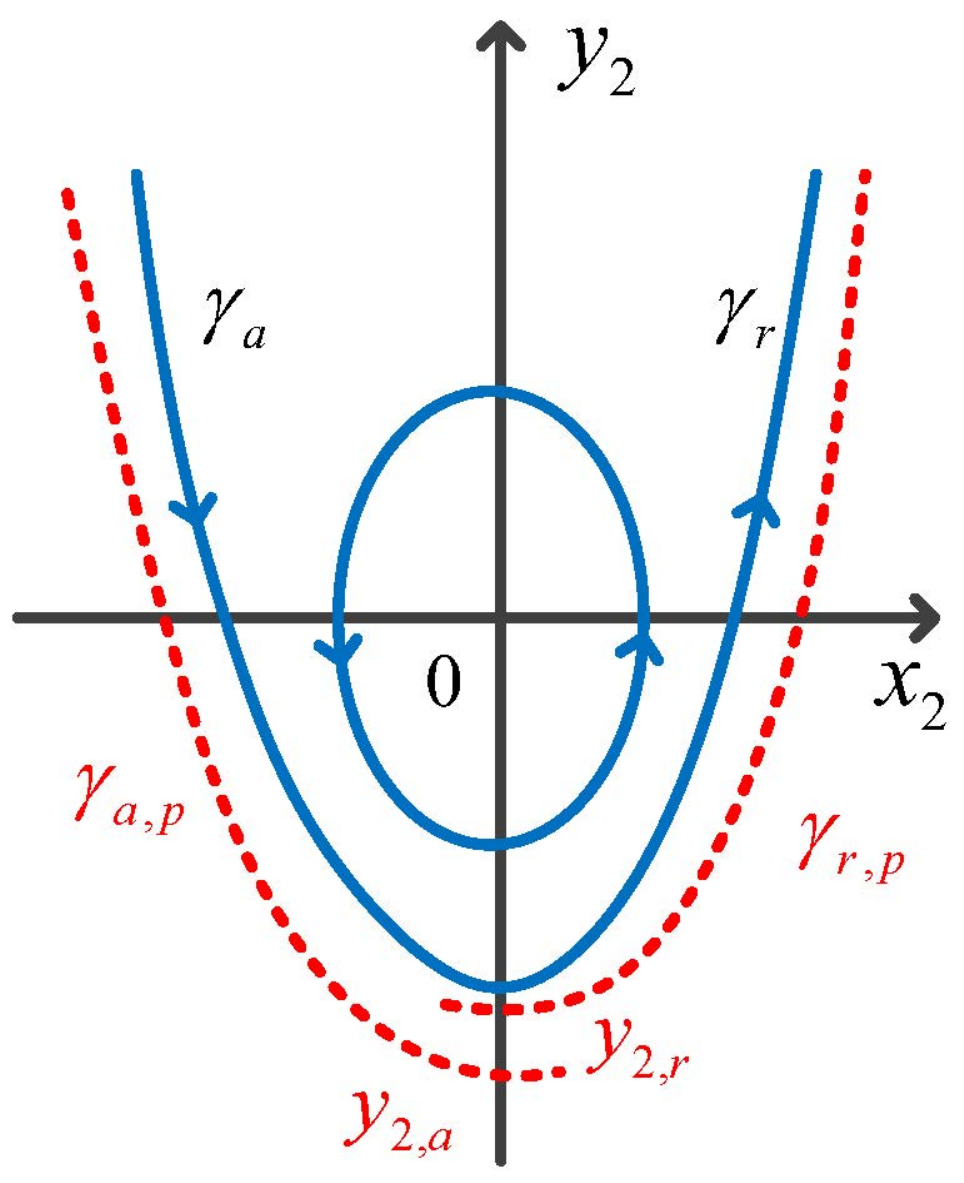}
\end{minipage}
\label{fg-level-curve}
}%
\centering
\caption{
\ref{fg-poincare-sphere} Phase portrait of system (\ref{int-1}) on the Poincar\'e disc.
\ref{fg-level-curve} The level curves of  $H$ and the perturbation of $\gamma$.
}
\end{figure}

In this section,
we establish the existence of canards near canard points
based on the blow-up technique and the Melnikov theory.
In order to study double canards near canard points $(\alpha_{j},\omega_{j})$,
we take a quasihomogeneous blow-up transformation $\Pi$ in the form
\begin{eqnarray*}
x=\bar{r}\bar{x}, \ \ y=\bar{r}^{2}\bar{y}, \ \ \varepsilon=\bar{r}^{2}\bar{\varepsilon},\ \
\lambda=\bar{r}\bar{\lambda}, \ \ \eta_{i}=\bar{r}\bar{\eta}_{i}, \ \ i=1,...,m.
\end{eqnarray*}
In the chart $K_{2}$ this blow-up transformation is reduced to $\Pi_{2}$ of the form
\begin{eqnarray}
\label{blow-up-2}
x=r_{2}x_{2}, \ \ y=r^{2}_{2}y_{2}, \ \ \varepsilon=r^{2}_{2},\ \
\lambda=r_{2}\lambda_{2}, \ \ \eta_{i}=r_{2}\eta_{2,i}, \ \ i=1,...,m.
\end{eqnarray}
By substituting (\ref{blow-up-2}) into (\ref{fast-normal-1}) and taking a time rescaling,
system (\ref{fast-normal-1}) is changed into
\begin{eqnarray}
\label{chart-2}
X_{2}(x_{2},y_{2},r_{2},\lambda_{2},\eta_{2}):\ \
\begin{split}
 x_{2}'  &=-y_{2}+x_{2}^{2}+r_{2}\left(a_{1}x_{2}-a_{2}x_{2}y_{2}+a_{3}x_{2}^{3}\right)
     +O\left(\!r_{2}(r_{2}+\lambda_{2}+\sum_{i=1}^{m}\eta_{2,i})\!\right),
\\
y_{2}'   &=x_{2}-\lambda_{2}+r_{2}\left(a_{4}x_{2}^{2}+a_{5}y_{2}\right)+\sum_{i=1}^{m}a_{5+i}\eta_{2,i}
+O\left(\!r_{2}(r_{2}+\lambda_{2}+\sum_{i=1}^{m}\eta_{2,i})\!\right),
\end{split}
\end{eqnarray}
where the constants $a_{i}$ are given by
\begin{eqnarray*}
&&a_{1}=\frac{\partial \phi_{3}}{\partial x}(0), \ \ \
a_{2}=\frac{\partial \phi_{1}}{\partial x}(0),\ \ \
a_{3}=\frac{\partial \phi_{2}}{\partial x}(0),\\
&&a_{4}=\frac{\partial \phi_{4}}{\partial x}(0),\ \ \
a_{5}=\phi_{6}(0),\ \ \
a_{5+i}=\phi_{6+i}(0), \ \ i=1,...,m.
\end{eqnarray*}
When $r_{2}=0$, $\lambda_{2}=0$, $\eta_{2}=0$,
system (\ref{chart-2}) is reduced to an integral system
\begin{eqnarray}
\label{int-1}
X_{2}^{I}(x_{2},y_{2}):\ \
\begin{split}
 x_{2}'  &= -y_{2}+x_{2}^{2},
\\
y_{2}'   &= x_{2}.
\end{split}
\end{eqnarray}
The solutions of this integral system are determined by the level curves of the function $H$, which is given by
\begin{eqnarray*}
H(x_{2},y_{2}):=\frac{1}{2}e^{-2y_{2}}\left(y_{2}-x_{2}^{2}+\frac{1}{2}\right), \ \ \ (x_{2},y_{2})\in \mathbb{R}^{2}.
\end{eqnarray*}
Clearly,
system (\ref{int-1}) has a solution $\gamma(t)$, $t\in \mathbb{R}$,
in the form
\begin{eqnarray}
\label{df-gamma}
\gamma(t)=\left(\frac{1}{2}t, \frac{1}{4}t^{2}-\frac{1}{2}\right), \ \ \ t\in\mathbb{R}.
\end{eqnarray}
By applying the Poincar\'e compactification (see, for instance, \cite[section V.1, p.321]{ZZF-etal}),
we can obtain the phase portrait of this integral system in the Poincar\'e disc,
which is shown in the Figure \ref{fg-poincare-sphere}.
We see that $\gamma$ is a heteroclinic orbit connecting two infinite equilibria.
Let  $\gamma_{a}=\{\gamma(t):t\leq 0\}$ and $\gamma_{r}=\{\gamma(t):t\geq 0\}$.
Consider the perturbed system (\ref{chart-2}) of the integral system (\ref{int-1}).
Assume that the heteroclinic orbit $\gamma$ is broken into $\gamma_{a,p}$ and $\gamma_{r,p}$.
By the continuous dependency on parameters,
the orbits $\gamma_{a,p}$ and $\gamma_{r,p}$ transversally intersects $y_{2}$-axis
at $(0,y_{2,a})$ and $(0,y_{2,r})$, which are in a small neighborhood of $(0,-1/2)$.
See Figure \ref{fg-level-curve}.
Clearly, the constants $y_{2,a}$ and $y_{2,r}$ depend on the parameters $r_{2},\lambda_{2}$ and $\eta_{2}$.
Here, we write $y_{2,a}$ and $y_{2,r}$  for simplicity.
To investigate the persistence of the heteroclinic orbit $\gamma$,
we define the so-called {\it distance function} $\mathcal{D}(r_{2},\lambda_{2},\eta_{2,1},...,\eta_{2,m})$ by
\begin{eqnarray*}
\mathcal{D}(r_{2},\lambda_{2},\eta_{2,1},...,\eta_{2,m})=y_{2,a}-y_{2,r}.
\end{eqnarray*}
If $\mathcal{D}(r_{2},\lambda_{2},\eta_{2,1},...,\eta_{2,m})=0$ for some suitable parameters,
then the heteroclinic orbit $\gamma$ is persistent.
Note that the heteroclinic orbit $\gamma$ in the integral system (\ref{int-1}) is unbounded,
then the major obstacle is to make sure
whether the distance function  $\mathcal{D}(r_{2},\lambda_{2},\eta_{2,1},...,\eta_{2,m})$
can be given by the classical Melnikov computation (see, for instance, \cite{Chow-Hale,Guck-Holmes-83,Melnikov}).
This problem can be solved by the method obtained by  Wechselberger in \cite{Wech-02}.
Roughly speaking,
for the extended system of (\ref{chart-2}), which is in the form
\begin{eqnarray}
\label{chart-sys-1}
\begin{split}
 x_{2}'  &=-y_{2}+x_{2}^{2}+r_{2}\left(a_{1}x_{2}-a_{2}x_{2}y_{2}+a_{3}x_{2}^{3}\right)
     +O\left(\!r_{2}(r_{2}+\lambda_{2}+\sum_{i=1}^{m}\eta_{2,i})\!\right),
\\
y_{2}'   &=x_{2}-\lambda_{2}+r_{2}\left(a_{4}x_{2}^{2}+a_{5}y_{2}\right)+\sum_{i=1}^{m}a_{5+i}\eta_{2,i}
+O\left(\!r_{2}(r_{2}+\lambda_{2}+\sum_{i=1}^{m}\eta_{2,i})\!\right),
\\
r_{2}'   &=0, \ \ \ \lambda_{2}'   =0, \ \ \ \eta_{2,i}'   =0, \ \ i=1,...,m,
\end{split}
\end{eqnarray}
if all solutions of the extended system (\ref{chart-sys-1}) near the heteroclinic orbit $\gamma$
are of at most algebraic growth for $t\to \pm \infty$.
Then  the  {\it distance function} $\mathcal{D}(r_{2},\lambda_{2},\eta_{2,1},...,\eta_{2,m})$
can be similarly obtained as in the classical case.
More precisely,
we have the following lemma.

\begin{figure}[!htp]
  \centering
  \includegraphics[width=5.8cm]{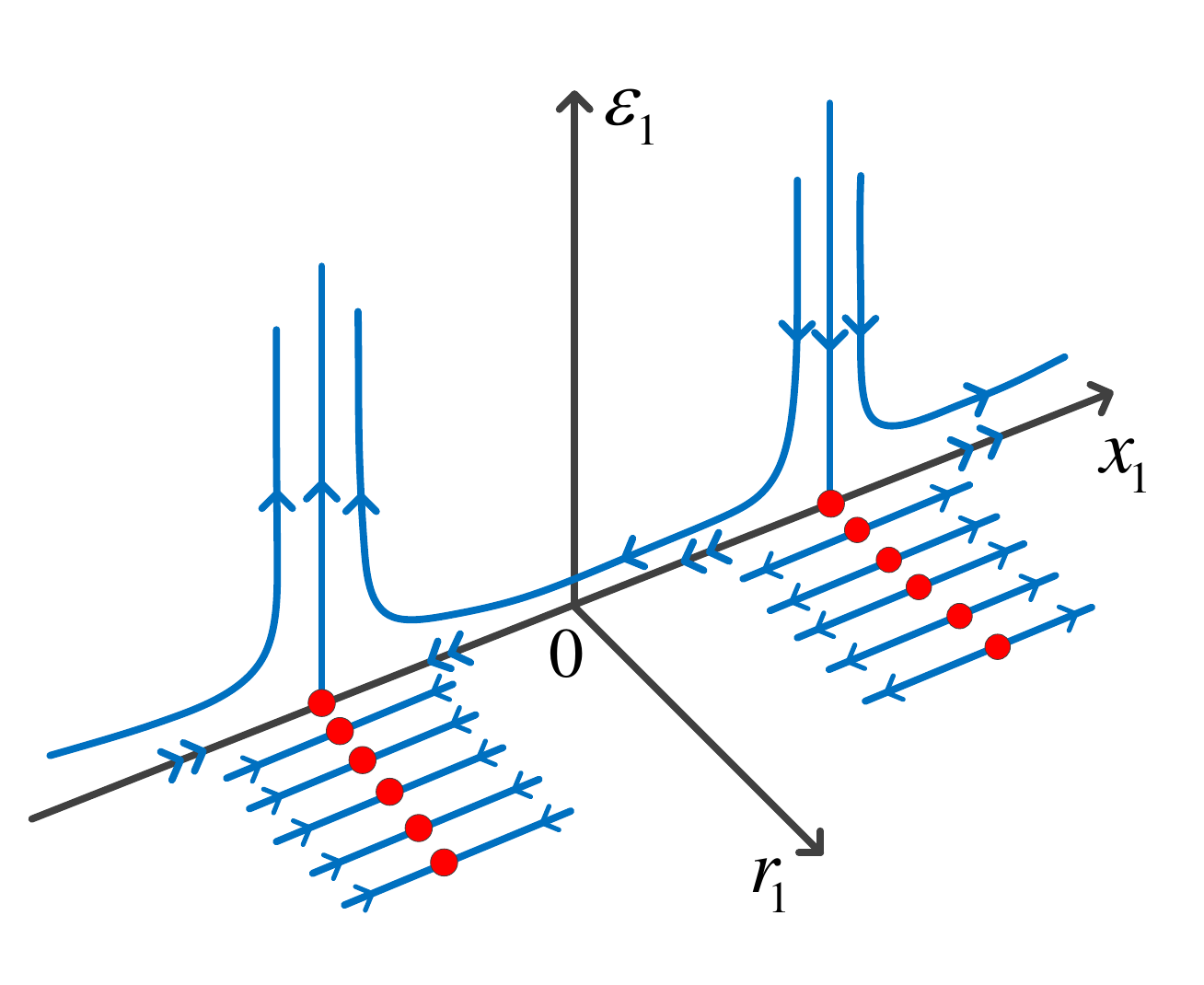}
  \caption{Dynamics of system (\ref{chart-sys-2}) with $\lambda_{1}=0$ and $\eta_{1,i}=0$ for $i=1,...,m$.}
  \label{fg-chart-1}
\end{figure}

\begin{lemma}
\label{lm-dist-funct}
The distance function $\mathcal{D}(r_{2},\lambda_{2},\eta_{2,1},...,\eta_{2,m})$ has the expansion
\begin{eqnarray*}
\mathcal{D}(r_{2},\lambda_{2},\eta_{2,1},...,\eta_{2,m})
=d_{r_{2}}r_{2}+d_{\lambda_{2}}\lambda_{2}+\sum_{i=1}^{m}d_{\eta_{2,i}}\eta_{2,i}+O\left(|(r_{2},\lambda_{2},\eta_{2,1},...,\eta_{2,m})|^{2}\right),
\end{eqnarray*}
where
\begin{eqnarray*}
d_{r_{2}}=-\frac{\sqrt{2\pi}}{8}(4a_{1}-a_{2}+3a_{3}-2a_{4}+2a_{5}),\ \ \
d_{\lambda_{2}}=-\sqrt{2\pi},  \ \ \ d_{\eta_{2,i}}=\sqrt{2\pi}a_{5+i},\ \ i=1,...,m.
\end{eqnarray*}
Furthermore,
the distance function $\mathcal{D}$ is $C^{k}$ smooth.
\end{lemma}
{\bf Proof.}
To apply the Melnikov theory in \cite[Theorem 1]{Wech-02},
we need to study the dynamics of the extended system (\ref{chart-sys-1}) at infinity.
By the transformation $\Pi_{12}$ in the form
\begin{eqnarray*}
\label{blow-up-1}
x_{2}=x_{1}\varepsilon_{1}^{-\frac{1}{2}}, \ \
y_{2}=\varepsilon_{1}^{-1},\ \
r_{2}=r_{1}\varepsilon_{1}^{\frac{1}{2}}, \ \
\lambda_{2}=\lambda_{1}\varepsilon_{1}^{-\frac{1}{2}},\ \
\eta_{2,i}=\eta_{1,i}\varepsilon_{1}^{-\frac{1}{2}}, \ \ i=1,...,m,\ \ \varepsilon_{1}>0,
\end{eqnarray*}
then this transformation give the blow-up transformation $\Pi$ in the chart $K_{1}$,
that is,
\begin{eqnarray}
x=r_{1}x_{1}, \ \ y=r^{2}_{1}, \ \ \varepsilon=r^{2}_{1}\varepsilon_{1},\ \
\lambda=r_{1}\lambda_{1}, \ \ \eta_{i}=r_{1}\eta_{1,i}, \ \ i=1,...,m.
\end{eqnarray}
By substituting (\ref{blow-up-1}) into (\ref{chart-sys-1})
and taking a rescaling of time $t\to \varepsilon_{1}^{1/2}t$,
the extended system (\ref{chart-sys-1}) is changed into
\begin{eqnarray}
\label{chart-sys-2}
\begin{split}
 x_{1}'  &=-1+x_{1}^{2}+r_{1}(a_{1}x_{1}\varepsilon_{1}-a_{2}x_{1}+a_{3}x_{1}^{3})
     -\frac{1}{2}\varepsilon_{1}x_{1}Q(x_{1},r_{1},\varepsilon_{1},\lambda_{1},\eta_{1})+{\rm h.o.t},
\\
r_{1}'   &=\frac{1}{2}\varepsilon_{1}r_{1}Q(x_{1},r_{1},\varepsilon_{1},\lambda_{1},\eta_{1}),
\\
\varepsilon_{1}'   &=-\varepsilon_{1}^{2}Q(x_{1},r_{1},\varepsilon_{1},\lambda_{1},\eta_{1}),
\\
\lambda_{1}'   &=-\frac{1}{2}\lambda_{1}\varepsilon_{1}Q(x_{1},r_{1},\varepsilon_{1},\lambda_{1},\eta_{1}),
\\
\eta_{1,i}'   &=-\frac{1}{2}\eta_{1,i}\varepsilon_{1}Q(x_{1},r_{1},\varepsilon_{1},\lambda_{1},\eta_{1}),
\ \ \ i=1,...,m,
\end{split}
\end{eqnarray}
where
\begin{eqnarray*}
Q(x_{1},r_{1},\varepsilon_{1},\lambda_{1},\eta_{1})
=x_{1}-\lambda_{1}+r_{1}(a_{4}x_{1}^{2}+a_{5})+{\rm  h.o.t},\ \ \
{\rm  h.o.t}=O\left(r_{1}(\lambda_{1}+\sum_{i=1}^{m}\eta_{1,i})\right).
\end{eqnarray*}
By a direct computation,
the points $E_{a}:=(-1,0,...,0)$ and $E_{r}:=(1,0,...,0)$ are two equilibria of the transformed system (\ref{chart-sys-2}),
which  have $(3+m)$-dimensional center manifolds $\mathcal{M}_{a}^{c}$ and $\mathcal{M}_{r}^{c}$,
and $1$-dimensional stable manifold
and $1$-dimensional unstable manifold along $x_{1}$-direction, respectively.  See Figure \ref{fg-chart-1}.
In the invariant plane $\{(x_1,0,\varepsilon_{1},0,..,0)\}$,
the transformed system (\ref{chart-sys-2}) is reduced to
\begin{eqnarray}
\label{reduce-2}
\begin{split}
x_{1}' &=-1+x_{1}^{2}-\frac{1}{2}\varepsilon_{1}x_{1}^{2},\\
\varepsilon_{1}' &=-\varepsilon_{1}^{2}x_{1}.
\end{split}
\end{eqnarray}
By \cite[Theorem 7.1, p.114]{ZZF-etal},
we obtain that the equilibrium $(-1,0)$ of system (\ref{reduce-2}) is a saddle-node with
an unique repelling center manifold $\tilde{\mathcal{M}}^{c}_{a}$ in $\varepsilon_{1}>0$,
and the equilibrium $(1,0)$ a saddle-node with
an unique attracting center manifold $\tilde{\mathcal{M}}^{c}_{r}$ in $\varepsilon_{1}>0$.
Note that $\Pi_{12}(\gamma_{a})=\tilde{\mathcal{M}}^{c}_{a}$ and $\Pi_{12}(\gamma_{r})=\tilde{\mathcal{M}}^{c}_{r}$,
then all solutions of the extended system (\ref{chart-sys-1}) near the algebraic growth $\gamma$
are of at most algebraic growth for $t\to \pm \infty$.

Recall that the heteroclinic orbit $\gamma$  has the form (\ref{df-gamma}),
then the  variational equation along $\gamma$ is in the form
$u'=tu-v$ and $v'=u$,
whose adjoint equation (see \cite[p.80]{Hale-80}) is given by $u'= -tu-v$ and $v'=u$.
Note that the collection of the solutions of this adjoint equation decaying
exponentially in forward and backward time is one-dimensional space,
which is spanned by $\tilde{\gamma}$ of the form
\begin{eqnarray*}
\tilde{\gamma}(t)=\left(-te^{-t^{2}/2}, e^{-t^{2}/2}\right), \ \ \ t\in\mathbb{R}.
\end{eqnarray*}
and the vector $\tilde{\gamma}(0)=(0,1)$,
then  by applying Theorem 1 and Proposition 1 in \cite{Wech-02},
the distance function $\mathcal{D}(r_{2},\lambda_{2},\eta_{2,1},...,\eta_{2,m})$ is $C^{k}$
and has the expansion
\begin{eqnarray*}
\mathcal{D}(r_{2},\lambda_{2},\eta_{2,1},...,\eta_{2,m})
=d_{r_{2}}r_{2}+d_{\lambda_{2}}\lambda_{2}+\sum_{i=1}^{m}d_{\eta_{2,i}}\eta_{2,i}+O\left(|(r_{2},\lambda_{2},\eta_{2,1},...,\eta_{2,m})|^{2}\right),
\end{eqnarray*}
where the coefficients $d_{r_{2}}$, $d_{\lambda_{2}}$ and $d_{\eta_{2,i}}$ are given by
\begin{eqnarray*}
d_{r_{2}}\!\!\!&=&\!\!\!\frac{\partial \mathcal{D}}{\partial r_{2}}(0)
=\int_{-\infty}^{+\infty}\langle\tilde{\gamma}(t), \frac{\partial X_{2}}{\partial r_{2}}(\gamma(t),0)\rangle\,dt
=-\frac{\sqrt{2\pi}}{8}(4a_{1}-a_{2}+3a_{3}-2a_{4}+2a_{5}),\\
d_{\lambda_{2}}\!\!\!&=&\!\!\!\frac{\partial \mathcal{D}}{\partial \lambda_{2}}(0)
=\int_{-\infty}^{+\infty}\langle\tilde{\gamma}(t),\frac{\partial X_{2}}{\partial \lambda_{2}}(\gamma(t),0)\rangle\,dt=-\sqrt{2\pi},\\
d_{\eta_{2,i}}\!\!\!&=&\!\!\!\frac{\partial \mathcal{D}}{\partial \eta_{2,i}}(0)
=\int_{-\infty}^{+\infty}\langle\tilde{\gamma}(t),\frac{\partial X_{2}}{\partial \eta_{2,i}}(\gamma(t),0)\rangle\,dt
=\sqrt{2\pi}a_{5+i},\ \ \ i=1,...,m,
\end{eqnarray*}
the vector field $X_{2}$ is  given by (\ref{chart-sys-1}),
the symbol $\langle\cdot,\,\cdot\rangle$ denotes the inner product of two vectors.
Therefore, the proof is now complete.
\hfill$\Box$

By Lemma \ref{lm-dist-funct},
we can obtain  the simultaneous occurrence of two canards near $(\alpha_{j},\omega_{j})$, $j=1,2$.

\begin{lemma}
\label{cor-canard-solu}
Suppose that system (\ref{slow-1}) satisfies the hypotheses {\bf (H1)}-{\bf (H5)},
and
\begin{eqnarray}
\label{con-rank-2}
{\rm Rank}
\left(
\begin{array}{llcl}
-f_{xx}(P_{1})(2g_{x}(P_{1}))^{-1}G_{1}d_{1,\lambda_{2}} & d_{1,\eta_{2,1}} & \cdot\cdot\cdot & d_{1,\eta_{2,m}}\\
-f_{xx}(P_{2})(2g_{x}(P_{2}))^{-1}G_{2}d_{2,\lambda_{2}} & d_{2,\eta_{2,1}}  & \cdot\cdot\cdot & d_{2,\eta_{2,m}}
\end{array}
\right)
=2.
\end{eqnarray}
Then there exists a $C^{k}$ function $\mu(\varepsilon)$
defined on $(0,\varepsilon_{1})$ for a small $\varepsilon_{1}>0$
such that two canards occur concurrently near the canard points $(\alpha_{j},\omega_{j})$, $j=1,2$.
\end{lemma}
{\bf Proof.}
 By Lemma \ref{lm-normal-form} we obtain that by the translations $\mathcal{T}_{j}$ of the form
\begin{eqnarray*}
\mathcal{T}_{j}(x,y,\lambda,\eta,\varepsilon)
      =\left(x+\tilde{x}_{j}((\lambda+\tilde{\lambda}_{j}(\varepsilon),\eta)+\mu_{0}),
      y+\tilde{y}_{j}((\lambda+\tilde{\lambda}_{j}(\varepsilon),\eta)+\mu_{0}),
   (\lambda+\tilde{\lambda}_{j}(\varepsilon),\eta)+\mu_{0},\varepsilon
      \right),
\end{eqnarray*}
and then by the rescaling $\mathcal{S}_{j}$ of the form
\begin{eqnarray*}
\mathcal{S}_{j}(x,y,\lambda,\eta,\varepsilon)
      =\left(\frac{2}{f_{xx}(P_{j})}x,\,
              -\frac{2}{f_{xx}(P_{j})f_{y}(P_{j})}y,\,
              -\frac{2g_{x}(P_{j})}{f_{xx}(P_{j})g_{\lambda}(P_{j})}\lambda,\,
              \eta,\,
              -\frac{1}{f_{y}(P_{j})g_{x}(P_{j})}\varepsilon
      \right),
\end{eqnarray*}
where $P_{j}=(\alpha_{j},\omega_{j},\mu_{0},0)$, and $\tilde{x}_{j}$, $\tilde{y}_{j}$ and $\tilde{\lambda}_{j}$  are respectively given by
(\ref{x-expan}), (\ref{y-expan}) and (\ref{lambda-expan}),
system (\ref{fast-1}) near the canard points $(\alpha_{j},\omega_{j})$ can be
changed into (\ref{fast-normal-1}).
Similarly to Lemma \ref{lm-dist-funct},
the distance functions $\mathcal{D}_{j}$ associated with the canard points $(\alpha_{j},\omega_{j})$
have the following expansions
\begin{eqnarray}
\label{expan-dis-j}
\mathcal{D}_{j}(r_{2},\lambda_{2},\eta_{2,1},...,\eta_{2,m})
=d_{j,r_{2}}r_{2}+d_{j,\lambda_{2}}\lambda_{2}+\sum_{i=1}^{m}d_{j,\eta_{i}}\eta_{2,i}+O\left(|(r_{2},\lambda_{2},\eta_{2,1},...,\eta_{2,m})|^{2}\right),
\end{eqnarray}
where  $d_{j,r_{2}}$, $d_{j,\lambda_{2}}$ and $d_{j,\eta_{2,i}}$ are given by (\ref{df-coeff}).
Let a projection $\mathcal{P}_{2}$
be defined by $\mathcal{P}_{2}(x_{2},y_{2},r_{2},\lambda_{2},\eta_{2})=(r_{2},\lambda_{2},\eta_{2})$.
Then to finish the proof, it suffices to solve the equations
$\mathcal{D}_{j}\circ\mathcal{P}_{2}\circ\Pi_{2}\circ\mathcal{S}_{j}\circ\mathcal{T}_{j}(x,y,\lambda,\eta,\varepsilon)=0$,
which are equivalent to the existence of solutions for the following  equations
\begin{eqnarray}
\label{eq-text}
\begin{split}
&-f_{y}(P_{j},\mu_{0},0)g_{x}(P_{j},\mu_{0},0)d_{j,r_{2}}\varepsilon
-\frac{f_{xx}(P_{j},\mu_{0},0)G_{j}}{2g_{x}(P_{j})}d_{j,\lambda_{2}}\left(\lambda-\lambda_{0}-\tilde{\lambda}_{j}(\varepsilon)\right)\\
&\ \  +\langle (d_{j,\eta_{2,1}},...,d_{j,\eta_{2,m}}),\eta-\eta_{0}\rangle
+o(|(\varepsilon,\mu-\mu_{0})|)=0, \ \ \ j=1,2,
\end{split}
\end{eqnarray}
where $G_{j}$ are defined as in {\bf (H4)}.
Since (\ref{con-rank-2}) holds,
then there exists a certain $i\in\{1,...,m\}$ such that
\begin{eqnarray*}
{\rm det}
\left(
\begin{array}{llcl}
-f_{xx}(P_{1})(2g_{x}(P_{1}))^{-1}G_{1}d_{1,\lambda_{2}} & d_{1,\eta_{2,i}}\\
-f_{xx}(P_{2})(2g_{x}(P_{2}))^{-1}G_{2}d_{2,\lambda_{2}} & d_{2,\eta_{2,i}}
\end{array}
\right)
\neq 0.
\end{eqnarray*}
This together with the {\it Implicit Function Theorem} yields that
there is an open interval $(0,\varepsilon_{1})$ for a small $\varepsilon_{1}>0$
and exactly two $C^{k}$ functions $\tilde{\lambda}(\varepsilon)$ and $\tilde{\eta}_{i}(\varepsilon)$
having the same expansions as in (\ref{df-lambda-eta}),
such that for each $\varepsilon\in (0,\varepsilon_{1})$ equations (\ref{eq-text}) have an solution
$
(\varepsilon,\lambda,\eta_{1},...,\eta_{m})
=(\varepsilon,\tilde{\lambda}(\varepsilon),\eta_{1,0},...,\eta_{i-1,0},\tilde{\eta}_{i}(\varepsilon),\eta_{i+1,0},...,\eta_{m,0}).
$
Thus, the proof is now complete.
\qquad $\Box$

\section{Proof of Theorem \ref{thm-exist}}
\label{sec-proof}
\setcounter{equation}{0}
\setcounter{lemma}{0}
\setcounter{theorem}{0}
\setcounter{remark}{0}

In this section we give the proof for Theorem \ref{thm-exist}
by Lemma \ref{lm-dist-funct} and \cite[Theorem 12]{Mamouhdi-Roussarie-12} (or \cite[Proposition 2]{Dumortier-Roussarie-07}).
We will see that there exists a codimension 2 limit cycle bifurcating from a double canard slow-fast cycle.
We refer the readers to \cite{Dumortier-Roussarie-08,Mamouhdi-Roussarie-12} for the precise definition of the codimension of a limit cycle.

{\bf Proof of Theorem \ref{thm-exist}.}
The proof for this theorem is divided into two steps.

{\bf Step 1.}
We first prove that under some suitable conditions,
there exists  a codimension 2 limit cycle bifurcating from $\Gamma(s_{1}^{*},s_{2}^{*})$,
from which three limit cycles can arise.
Since (\ref{con-rank-2}) holds,
then there exists a local diffeomorphism $\mathcal{F}$
transforming $\vartheta:=(\lambda,\eta_{i})$ near $(\lambda_{0},\eta_{i,0})$
to $(\beta_{1},\beta_{2})=(\beta_{1}(\vartheta),\beta_{2}(\vartheta))$ such that
the canard point $(\alpha_{1},\omega_{1})$ appears when $\beta_{1}=0$,
and the other  canard point $(\alpha_{2},\omega_{2})$ appears when $\beta_{2}=0$.
Since the curves
$\left\{(s_{1},s_{2})\in\Omega: \mathcal{I}_{1}(s_{1}) +\mathcal{I}_{2}(s_{2})=0 \right\}$
and $\left\{(s_{1},s_{2})\in\Omega: \mathcal{I}_{3}(s_{2}) +\mathcal{I}_{4}(s_{1})=0 \right\}$
transversally intersects at the point $(s_{1}^{*},s_{2}^{*})\in\Omega$,
then the double canard slow-fast cycle $\Gamma(s_{1}^{*},s_{2}^{*})$ is codimension 2 (see \cite[Definition 7]{Mamouhdi-Roussarie-12}).
By \cite[Theorem 12]{Mamouhdi-Roussarie-12}
there exists a sufficiently small $\varepsilon_{0}$ and  two continuous functions $\beta_{1}(\varepsilon)$ and $\eta_{i}(\varepsilon)$ for
$0\leq \varepsilon<\varepsilon_{0}$ such that
system (\ref{fast-1}) with $\beta_{1}=\beta_{1}(\varepsilon)$, $\eta_{i}=\eta_{i}(\varepsilon)$
and  fixed $\eta_{j}=\eta_{j,0}$ for $j\neq i$
has a codimension 2 limit cycle $\Gamma_{\varepsilon}(s_{1}^{*},s_{1}^{*})$,
which can produce three hyperbolic limit cycles for some suitable parameters.

\begin{figure}[!htp]
  \centering
  \includegraphics[width=7.5cm]{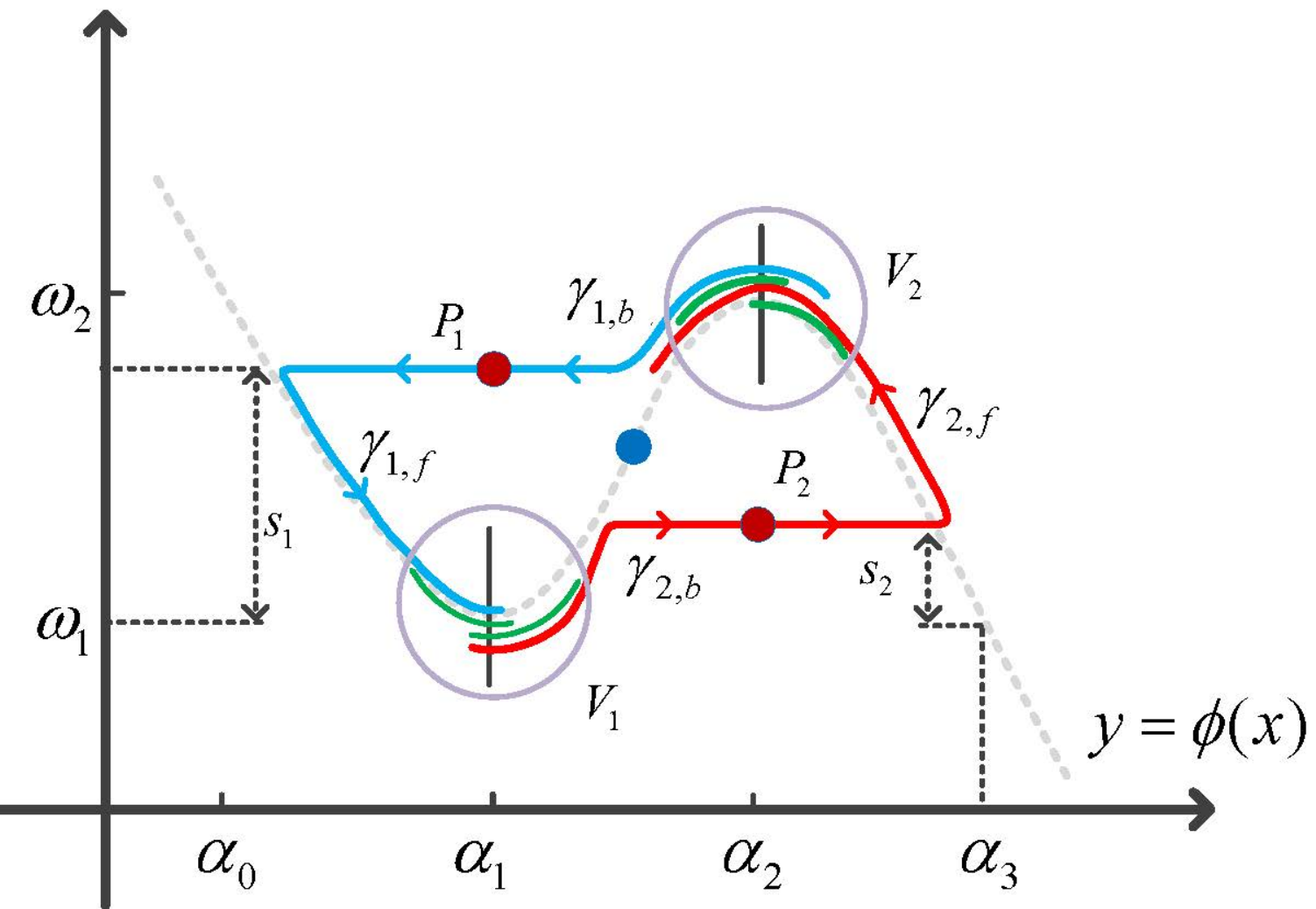}
  \caption{Perturbation of double canard slow-fast cycle $\Gamma(s_{1},s_{2})$.}
  \label{fg-butterfly-pert}
\end{figure}

{\bf Step 2.}
Secondly,
we give the explicit expansions of $\lambda(\varepsilon)$ and $\eta_{i}(\varepsilon)$
with  $(\lambda(\varepsilon),\eta_{i}(\varepsilon))=\mathcal{F}^{-1}((\beta_{1}(\varepsilon),\eta_{i}(\varepsilon)))$
for $0\leq \varepsilon<\varepsilon_{0}$.
Fix $\eta_{j}=\eta_{j,0}$ for $j\neq i$ and set $P_{j}=(\alpha_{j},\omega_{1}+s_{j}^{*})$, $j=1,2$.
Let the forward orbits and the backward orbits of $P_{j}$ under the flow of system (\ref{fast-1})
be respectively denoted by $\gamma_{j,f}(t), t\geq 0$,  and $\gamma_{j,b}(t)$, $t\leq 0$,
which satisfy $\gamma_{j,f}(0)=\gamma_{j,b}(0)=P_{j}$.
For each $j=1,2$,
we take a small open neighborhood $V_{j}$ of the canard points $(\alpha_{j},\omega_{j})$
such that near the point $(\alpha,\omega)=(\alpha_{j},\omega_{j})$,
system (\ref{fast-1}) can be changed into (\ref{fast-normal-1}).
By the Fenichel Theorem \cite[Theorem 9.1]{Fenichel-79},
there exists two open sets $I_{j,f}:=(a_{j,f}^{+},b_{j,f}^{+})\subset\mathbb{R}^{+}$
and $I_{j,b}:=(a_{j,b}^{-},b_{j,b}^{-})\subset\mathbb{R}^{-}$ such that
$\gamma_{j,f}(t)\in V_{j}$ for $t\in I_{j,f}$ and $\gamma_{j,f}(t)\notin V_{j}$ for $0\leq t\leq a_{j,f}^{+}$,
$\gamma_{1,b}(t)\in V_{2}$ for $t\in I_{1,b}$ and $\gamma_{1,b}(t)\notin V_{2}$ for $b_{1,b}^{-}\leq t\leq 0$,
and $\gamma_{2,b}(t)\in V_{1}$ for $t\in I_{2,b}$ and $\gamma_{2,b}(t)\notin V_{1}$ for $b_{2,b}^{-}\leq t\leq 0$.
See Figure \ref{fg-butterfly-pert}.

Consider system (\ref{fast-1}) in the sets $V_{j}$.
Following the  discussions as in the section \ref{sec-canard solu},
 system (\ref{fast-1}) can be transformed into
the similar systems as (\ref{chart-sys-1}).
Let these two transformed systems  be respectively denoted by $X_{2,j}(x_{2},y_{2},r_{2},\lambda_{2},\eta_{2})$,
which are the perturbations of the integral system $X^{I}_{2}(x_{2},y_{2})$.
To distinguish two canard points $(\alpha_{j},\omega_{j})$, $j=1,2$,
let the heteroclinic orbit $\gamma=\gamma_{a}\cup\gamma_{r}$ corresponding to $(\alpha_{j},\omega_{j})$
be denoted by $\gamma^{j}=\gamma_{a}^{j}\cup\gamma_{r}^{j}$.
Define the perturbations of $\gamma_{a}^{j}$ and $\gamma_{r}^{j}$ by
$\gamma_{a,p}^{j}$ and $\gamma_{r,p}^{j}$,
which respectively transversally intersect $x_{2}=0$ at
$(0,y_{2,a}^{j})$ and $(0,y_{2,r}^{j})$.
Let  the intersection point of  the transformed orbits of $\gamma_{j,f}$ (resp. $\gamma_{j,b}$) and $x_{2}=0$
be denoted by $(0,y_{f}^{j})$ (resp. $(0,y_{b}^{j})$).
By  \cite[Theorem 9.1]{Fenichel-79},
 the slow manifolds are exponentially close to normally hyperbolic branches of the critical manifold $\mathcal{C}_{0}$.
In the sets $V_{j}$, consider system (\ref{chart-sys-2}) in the chart $K_{1}$
and system (\ref{chart-sys-1}) in the chart $K_{2}$.
By the similar argument as in \cite[Lemma 5.1]{Krupa-Szmolyan-01SIMA},
we can obtain that there exists a positive constant $\kappa$ such that
\begin{eqnarray}
\label{est-1}
\begin{split}
&|y_{2,a}^{1}-y_{f}^{1}|=O(e^{-\kappa/r^{2}_{2}}), \ \
|y_{2,r}^{1}-y_{b}^{2}|=O(e^{-\kappa/r^{2}_{2}}),\\
&|y_{2,a}^{2}-y_{f}^{2}|=O(e^{-\kappa/r^{2}_{2}}),\ \
|y_{2,r}^{2}-y_{b}^{1}|=O(e^{-\kappa/r^{2}_{2}}),
\end{split}
\end{eqnarray}
and the partial derivatives of $y_{2,a}^{1}-y_{f}^{1}$, $y_{2,r}^{1}-y_{b}^{2}$,
$y_{2,a}^{2}-y_{f}^{2}$ and $y_{2,r}^{2}-y_{b}^{1}$ with respect to $r_{2}$, $\lambda_{2}$ and $\eta_{2,i}$, $i=1,...,m$,
have the similar estimates as above.
Note that  there exists a double canard cycle bifurcating from $\Gamma(s_{1}^{*},s_{2}^{*})$
if and only if
\begin{eqnarray*}
\tilde{\mathcal{D}}_{1}(r_{2},\lambda_{2},\eta_{2,i}):=y^{1}_{f}-y^{2}_{b}=0, \ \ \
\tilde{\mathcal{D}}_{2}(r_{2},\lambda_{2},\eta_{2,i}):=y^{2}_{f}-y^{1}_{b}=0.
\end{eqnarray*}
Set $\mathcal{D}_{j}(r_{2},\lambda_{2},\eta_{2,i}):=y_{2,a}^{j}-y_{2,r}^{j}$.
Note that
\begin{eqnarray*}
\tilde{\mathcal{D}}_{1}(r_{2},\lambda_{2},\eta_{2,i})
\!\!\!&=&\!\!\! \mathcal{D}_{1}(r_{2},\lambda_{2},\eta_{2,i})+(y_{f}^{1}-y_{2,a}^{1})+(y_{2,r}^{1}-y_{b}^{2}),\\
\tilde{\mathcal{D}}_{2}(r_{2},\lambda_{2},\eta_{2,i})
\!\!\!&=&\!\!\! \mathcal{D}_{2}(r_{2},\lambda_{2},\eta_{2,i})+(y_{f}^{2}-y_{2,a}^{2})+(y_{2,r}^{2}-y_{b}^{1}),
\end{eqnarray*}
then
similarly to the proof for Lemma \ref{cor-canard-solu},
there exists a double canard cycle bifurcating from $\Gamma(s_{1}^{*},s_{2}^{*})$
if and only if two equations with the expansions (\ref{eq-text}) has solutions.
Since (\ref{con-rk-2}) holds,
then letting $\eta_{j}=\eta_{j,0}$ for $j\neq i$,
by the {\it Implicit Function Theorem}
there exists a sufficiently small $\varepsilon_{0}$ and a unique continuous curve
$\mu(\varepsilon):=(\varepsilon,\tilde{\lambda}(\varepsilon),\eta_{1,0},...,\eta_{i-1,0},\tilde{\eta}_{i}(\varepsilon),\eta_{i+1,0},...,\eta_{m,0})$
for $0\leq \varepsilon<\varepsilon_{0}$
such that (\ref{eq-text}) has the solution $\mu=\mu(\varepsilon)$,
where $\tilde{\lambda}(\varepsilon)$ and $\tilde{\eta}_{i}(\varepsilon)$ have the same expansions defined as in  (\ref{df-lambda-eta}).
Then by the uniqueness we have that
$(\lambda(\varepsilon),\eta_{i}(\varepsilon))=\mathcal{F}^{-1}((\beta_{1}(\varepsilon),\eta_{i}(\varepsilon)))
=(\tilde{\lambda}(\varepsilon),\tilde{\eta}_{i}(\varepsilon))$ for $0\leq \varepsilon<\varepsilon_{0}$.
Thus, the proof is finished.
\hfill$\Box$

\section{Applications to Li\'{e}nard equations with quadratic damping}
\label{sec-lienard}
\setcounter{equation}{0}
\setcounter{lemma}{0}
\setcounter{theorem}{0}
\setcounter{remark}{0}

Consider a singularly perturbed Li\'{e}nard equation of the form
\begin{eqnarray}
\label{lienard-1}
\begin{split}
\frac{d x}{d t} &= x'= F(x)-y,
\\
\frac{d y}{d t}&= y'=  \varepsilon (\eta+\lambda x-F(x)),
\end{split}
\end{eqnarray}
where $(x,y)\in\mathbb{R}^{2}$, the function $F$ is a cubic polynomial with respect to $x$,
the parameters $\eta$, $\lambda$ and $\varepsilon$ are real and $\varepsilon>0$ is sufficiently small.
Li\'{e}nard equation (\ref{lienard-1}) is equivalent to the following equation
\begin{eqnarray*}
x''+F'(x)x'+P(x)=0,
\end{eqnarray*}
where the derivative $F'$ of the cubic polynomial $F$ and the function $P$ defined by $P(x):=-\varepsilon (\eta+\lambda x-F(x))$
are polynomials of degree two and three respectively.
Then system (\ref{lienard-1}) is always called a {\it cubic Li\'{e}nard equation with quadratic damping}
(see, for instance, \cite{Dumortier-Li-97,Dumortier-etal-00}).

We can compute that for each $(x,y)\in\mathbb{R}^{2}$,
\begin{eqnarray}\label{Bdixson-1}
\frac{\partial}{\partial x}(F(x)-y)+\frac{\partial}{\partial y}(\varepsilon (\eta+\lambda x-F(x)))=F'(x).
\end{eqnarray}
If the cubic polynomial $F$ satisfies that either $F'(x)\geq 0$ or $F'(x)\leq 0$ holds for all $x\in\mathbb{R}$,
then by {\it Bendixson's Theorem} (see, for instance, \cite[Theorem 7.10, p.188]{Dumortieretal-06}),
we obtain that (\ref{lienard-1}) has no limit cycles in $\mathbb{R}^{2}$.
Thus one essential assumption for the existence of limit cycles is that
the cubic polynomial $F$ has precisely two different extreme points.
This assumption implies that the critical manifold associated with singularly perturbed Li\'{e}nard equation  (\ref{lienard-1})
is $S$-shaped.
In this case, Li\'{e}nard equation (\ref{lienard-1}) can be normalized into a simpler form,
the results are summarized in the following lemma.
\begin{lemma}
\label{lm-simp}
Assume that the cubic polynomial $F$ in (\ref{lienard-1}) has precisely two different extreme points.
Then (\ref{lienard-1}) can be changed into
\begin{eqnarray}
\label{lienard-2}
\begin{split}
\frac{d x}{d t} &=x'= F(x)-y,
\\
\frac{d y}{d t}&=y'= \varepsilon \left(\eta+\lambda \left(x-\frac{1}{2}\right)-F(x)\right),
\end{split}
\end{eqnarray}
where the function $F$ satisfies the following:
\begin{eqnarray}
\label{df-F-0}
F(x)=-\frac{1}{3}x^3+\frac{1}{2}x^{2}, \ \ \ \
F'(x)=-x(x-1),\ \ \ \  x\in\mathbb{R}.
\end{eqnarray}
\end{lemma}
{\bf Proof.}
By  a translation
we can move the left extreme point to the origin.
Then we assume that the cubic polynomial $F$  satisfies that
\begin{eqnarray}
\label{df-F}
F(x)=\frac{1}{3}u x^{3}-\frac{1}{2}u\nu x^{2},  \ \ \ F'(x)=u x(x-\nu),
\end{eqnarray}
where the parameters $u$ and $\nu$ respectively satisfy $u\neq 0$ and $\nu>0$.
We see that  (\ref{lienard-1}) is equivalent to
\begin{eqnarray*}
\frac{d x}{d t} \!\!\!&=&\!\!\! z,
\\
\frac{d z}{d t} \!\!\!&=&\!\!\! -\varepsilon(\eta+\lambda x-F(x))-F'(x)z.
\end{eqnarray*}
Let $x=\nu \bar{x}$, $z=\nu \bar{z}$, $t=\bar{t}$,
$u=-\nu^{-2}\bar{u}$, $\nu=\bar{\nu}$, $\varepsilon=\bar{\varepsilon}$,
$\eta=\nu\bar{\eta}$ and $\lambda=\bar{\lambda}$.
After dropping the bars,
we can verify that (\ref{lienard-1}) is equivalent to
\begin{eqnarray*}
\frac{d x}{d t} \!\!\!&=&\!\!\! z,
\\
\frac{d z}{d t} \!\!\!&=&\!\!\! -\varepsilon\left(\eta+\lambda x+\frac{1}{3}u x^{3}-\frac{1}{2}u x^{2}\right)+u x(x-1)z.
\end{eqnarray*}
By taking the changes
$x=\widetilde{x}$, $z=u\widetilde{z}$, $t=\widetilde{t}/u$,
$\varepsilon=u\widetilde{\varepsilon}$, $\eta=u (\widetilde{\eta}-\widetilde{\lambda}/2)$,
$\lambda=u \widetilde{\lambda}$ and $u=\widetilde{u}$ in the above system,
and then dropping all tildes,
the parameter $u$ in the above system is changed to $u=1$.
Then  (\ref{lienard-1})  is equivalent to (\ref{lienard-2}) with $F$ given by (\ref{df-F}).
Thus, the proof is finished.
\hfill $\Box$

By a time rescaling $s=\varepsilon t$,
the corresponding  slow system of (\ref{lienard-2}) is in the form
\begin{eqnarray}
\label{lienard-slow-1}
\begin{split}
\varepsilon \frac{d x}{d s} &= \varepsilon \dot x=F(x)-y:=f(x,y),
\\
\frac{d y}{d s} &=\dot y=\left(\eta+\lambda \left(x-\frac{1}{2}\right)-F(x)\right):=g(x,\eta,\lambda).
\end{split}
\end{eqnarray}
The slow motion of system (\ref{lienard-2})
along the critical manifold
\begin{eqnarray*}
\mathcal{C}_{0}=\left\{(x,y)\in\mathbb{R}: y=-\frac{1}{3}x^3+\frac{1}{2}x^{2}\right\},
\end{eqnarray*}
is governed by
\begin{eqnarray}
\label{eq-slow-motion}
F'(x)\frac{d x}{ds}=g(x,\eta,\lambda).
\end{eqnarray}
At the points $(0,0)$ and $(1,1/6)$,
the function $f$ satisfies that
\begin{eqnarray*}
f(0,0)=f(1,1/6)=\frac{\partial f}{\partial x}(0,0)=\frac{\partial f}{\partial x}(1,1/6)=0,
\end{eqnarray*}
and the following nondegeneracy conditions:
\begin{eqnarray*}
\frac{\partial^{2} f}{\partial x^{2}}(1,1/6)=-1, \ \ \
\frac{\partial^{2} f}{\partial x^{2}}(0,0)=\frac{\partial f}{\partial y}(0,0)=\frac{\partial f}{\partial y}(1,1/6)=1.
\end{eqnarray*}
Clearly, there are two canard points on the critical manifold $\mathcal{C}_{0}$
if and only if the parameters $\eta$ and $\lambda$ respectively satisfy $\eta=1/12$ and $\lambda=1/6$.
In this case, at the points $(0,0)$ and $(1,1/6)$ we have that
\begin{eqnarray*}
&&\frac{\partial g}{\partial x}(0,1/12,1/6)=\frac{\partial g}{\partial x}(1,1/12,1/6)=1/6,\\
&&\frac{\partial g}{\partial \eta}(0,1/12,1/6)=\frac{\partial g}{\partial \eta}(1,1/12,1/6)=1,\\
&&\frac{\partial g}{\partial \lambda}(0,1/12,1/6)=-1/2, \ \
\frac{\partial g}{\partial \lambda}(1,1/12,1/6)=1/2.
\end{eqnarray*}
Then for $\eta=1/12$ and $\lambda=1/6$,
system (\ref{lienard-slow-1}) satisfies the hypotheses {\bf (H1)}-{\bf (H5)} stated in the section 2.

\begin{proposition}
\label{lm-lienard-butterfly}
Consider the Li\'enard system of the form (\ref{lienard-2}).
Then for sufficiently small $\varepsilon$,
there exists a function $\bar{\lambda}$ of the form
$\bar{\lambda}(\varepsilon)=-2\varepsilon/9+1/6+O(\varepsilon^{\frac{3}{2}})$
and $\eta=1/12$ such that double canards appears.
Furthermore, assume that system (\ref{lienard-2}) satisfies the conditions in Theorem \ref{thm-exist}.
Then there exist three big limit cycles which are hyperbolic in system (\ref{lienard-2}).
\end{proposition}
{\bf Proof.}
For simplicity,
let the pair $(x_{0},y_{0})$ denote either $(0,0)$ or  $(1,1/6)$.
By applying the changes
\begin{eqnarray}
x=2\zeta \tilde{x}+x_{0}, \ \ y=2\zeta \tilde{y}+y_{0}, \ \
\lambda=\frac{2}{3}\tilde{\lambda}+\frac{1}{6},\ \
\eta=\tilde{\eta}+\frac{1}{12}, \ \ \varepsilon=6\tilde{\varepsilon},
\label{app-tran-1}
\end{eqnarray}
where $\zeta=1$ for $(x_{0},y_{0})=(0,0)$ and $\zeta=-1$ for $(x_{0},y_{0})=(1,1/6)$,
we obtain the normal form stated as in Lemma \ref{lm-normal-form}, that is,
\begin{eqnarray}
\label{trans-lienard-1}
\begin{split}
\tilde{x}' &=-\tilde{y}+\tilde{x}^{2}\left(1-\frac{4}{3}\tilde{x}\right),
\\
\tilde{y}' &=\tilde{\varepsilon}\left(\tilde{x}(8\tilde{x}^{2}-6\tilde{x}+1)-\tilde{\lambda}(1-4\tilde{x})+3\zeta\tilde{\eta}\right).
\end{split}
\end{eqnarray}
In the chart $K_{2}$ with the blow-up transformation
\begin{eqnarray}
\tilde{x}=r_{2}x_{2}, \ \ \tilde{y}=r^{2}_{2}y_{2}, \ \ \tilde{\varepsilon}=r^{2}_{2},\ \
\tilde{\lambda}=r_{2}\lambda_{2}, \ \ \tilde{\eta}=r_{2}\eta_{2},
\label{app-tran-2}
\end{eqnarray}
by Lemma \ref{lm-dist-funct}
the distance functions $\mathcal{D}_{\pm}(r_{2},\lambda_{2},\eta_{2})$ associated with systems (\ref{trans-lienard-1})
have the expansions
\begin{eqnarray*}
\mathcal{D}_{\pm}(r_{2},\lambda_{2},\eta_{2})
=-2\sqrt{2\pi}r_{2}-\sqrt{2\pi}\lambda_{2}+3\sqrt{2\pi}\zeta\eta_{2}+O(|(r_{2},\lambda_{2},\eta_{2})|^{2}),
\end{eqnarray*}
where $\mathcal{D}_{+}$ and $\mathcal{D}_{-}$ correspond to the cases $\zeta=1$ and $\zeta=-1$, respectively.
To obtain the existence of  double canard cycles,
we consider the equations $\mathcal{D}_{\pm}(r_{2},\lambda_{2},\eta_{2})=0$.
Then by Lemma \ref{cor-canard-solu},
there exists an open neighbourhood $U_{r_{2}}(0)\subset\mathbb{R}$ of $r_{2}=0$
and exactly two smooth functions $\bar{\lambda}_{2}$ and $\bar{\eta}_{2}$ in the form
\begin{eqnarray*}
\bar{\lambda}_{2}(r_{2})=-2r_{2}+O(r_{2}^{2}), \ \ \ \bar{\eta}_{2}(r_{2})=O(r_{2}^{2}),
\end{eqnarray*}
such that
$\mathcal{D}_{\pm}(r_{2},\bar{\lambda}_{2}(r_{2}),\bar{\eta}_{2}(r_{2}))=0$ for each $r_{2}\in U_{r_{2}}(0)$.

We claim that $\bar{\eta}_{2}(r_{2})\equiv 0$ for each $r_{2}\in U_{r_{2}}(0)$.
Note that  the normal forms (\ref{trans-lienard-1}) of system (\ref{lienard-2}) near the points $(0,0)$ and $(1,1/6)$
have the only difference in the coefficient of $\bar{\eta}$.
Then by this symmetry we obtain that $\bar{\eta}_{2}(r_{2})=-\bar{\eta}_{2}(r_{2})$ for each $r_{2}\in U_{r_{2}}(0)$.
This implies that the claim holds.
Recalling the transformations (\ref{app-tran-1}) and (\ref{app-tran-2}),
and applying Theorem \ref{thm-exist},
we obtain this proposition.
Thus, the proof is finished.
\hfill $\Box$

\begin{figure}[!htp]
  \centering
  \includegraphics[width=8cm]{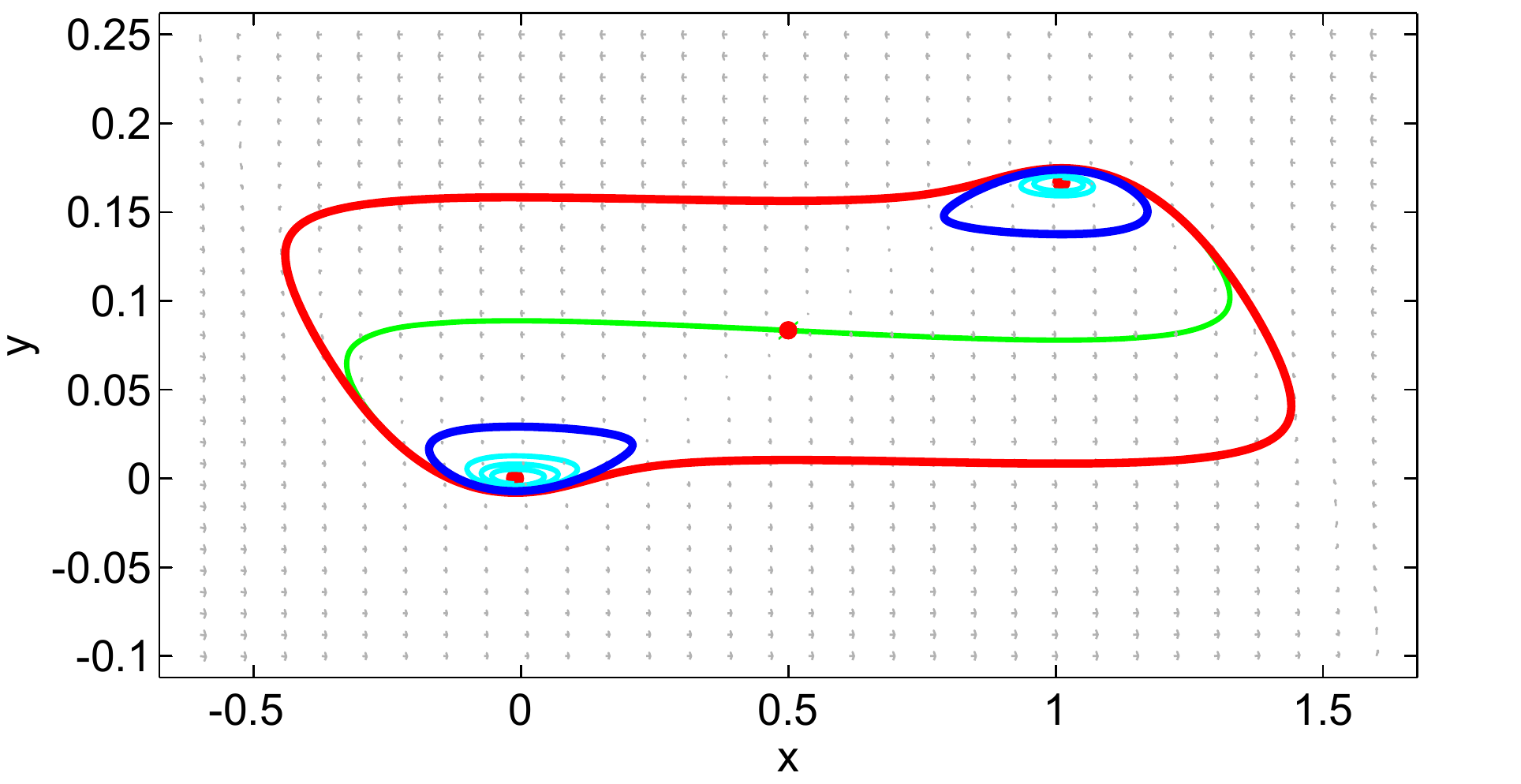}
  \caption{Coexistence of three limit cycles for system (\ref{lienard-2}) with $\varepsilon=0.045$, $\eta=1/12$ and $\lambda=0.163$.}
  \label{Fig-simu}
\end{figure}
\begin{remark}
\label{rmk-1}
It is highly difficult to analyze the intersections of the following curves
\begin{eqnarray*}
\left\{(s_{1},s_{2})\in\Omega: \mathcal{I}_{1}(s_{1}) +\mathcal{I}_{2}(s_{2})=0 \right\}, \ \ \
\left\{(s_{1},s_{2})\in\Omega: \mathcal{I}_{3}(s_{2}) +\mathcal{I}_{4}(s_{1})=0 \right\},
\end{eqnarray*}
where $\mathcal{I}_{i}$ are similarly defined as in section \ref{sec-mainresults}.
The detailed study of them will be given in the future work.
Numerical  simulation shows that there exists a big double canard cycle enclosing
two small limit cycles, each of which surrounds an stable focus. See Figure \ref{Fig-simu}.
Based on Proposition \ref{lm-lienard-butterfly} and Figure \ref{Fig-simu},
we conjecture that there exist three big limit cycles enclosing two small limit cycles under some suitable conditions.
If this conjecture is right,
then the maximal number of equation (\ref{lienard-1}) is at least five.
We hope that the results obtained in this paper are useful to explain this conjecture and
the  phenomenon in Figure \ref{Fig-simu}.
\end{remark}

\end{document}